\documentclass[3p]{elsarticle}   
\usepackage[bookmarks=false]{hyperref}
\usepackage{color} 
\usepackage{psfrag} 
\usepackage{subfig}
\usepackage{amsmath}
\usepackage{amsfonts}
\usepackage{bm}
\usepackage{graphics}
 \usepackage{pstool} 
\usepackage{todonotes}
\usepackage[]{algorithm2e}

\usepackage{pgfplots}
\usepackage{pgfplotstable}
\usepackage{filecontents}

\usepackage{tikz}
\usepackage{tikz-dimline}

\usetikzlibrary{spy}
\usetikzlibrary{decorations,arrows,shapes}    
\usetikzlibrary{backgrounds}  
\usetikzlibrary{positioning}
\usetikzlibrary{patterns}
\usetikzlibrary{calc} 

\usepackage{booktabs}   
\usepackage{makecell} 
\usepackage{colortbl}   

\usepackage{xspace}   
\usepackage{setspace}   
 
\usepackage{soul}
\usepackage{ulem}
\usepackage{currfile}
\usepackage{import}

\usepackage{cleveref}
\usepackage{dsfont}

\usepackage{amssymb}
\usepackage{grffile}
             
\usepackage{tikzscale}
\usepackage{tikz-3dplot}
\pgfplotsset{compat=1.8}
    
\usetikzlibrary{3d}
\usetikzlibrary{snakes}    

\usepackage{gensymb}
\usepackage{framed}

\usepackage[utf8]{inputenc} 
\usepackage[T1]{fontenc}

\pgfplotsset{compat=1.7}

\newcommand{\findmax}[3]{
    \pgfplotstablesort[sort key={#2},sort cmp={float >}]{\sorted}{#1}%
    \pgfplotstablegetelem{0}{#2}\of{\sorted}%
    \let #3=\pgfplotsretval%
}

\pgfplotscreateplotcyclelist{myCycleListBlackAndWhite}
{%
  {black, mark=o},
  {black, mark=square},
  {black, mark=triangle},
  {black, mark=diamond}, 
  {black, mark=pentagon}, 
  {black, mark=*},
  {black, mark=square*},
  {black, mark=triangle*},
  {black, mark=diamond*}, 
  {black, mark=pentagon*}, 
}

\definecolor{darkgreen}{rgb}{0,0.4,0} 
\definecolor{darkbrown}{rgb}{0.5, 0.396, 0.09}

\definecolor{c1}{rgb}{0.0, 0.4196078431372549, 0.6431372549019608}
\definecolor{c2}{rgb}{1.0, 0.5019607843137255, 0.054901960784313725}
\definecolor{c3}{rgb}{0.6705882352941176, 0.6705882352941176,
0.6705882352941176} \definecolor{c}{rgb}{0.34901960784313724, 0.34901960784313724, 0.34901960784313724}
\definecolor{c4}{rgb}{0.37254901960784315, 0.6196078431372549,
0.8196078431372549} \definecolor{c}{rgb}{0.7843137254901961, 0.3215686274509804, 0.0}
\definecolor{c5}{rgb}{0.5372549019607843, 0.5372549019607843,
0.5372549019607843} \definecolor{c}{rgb}{0.6352941176470588, 0.7843137254901961, 0.9254901960784314}
\definecolor{c6}{rgb}{1.0, 0.7372549019607844, 0.4745098039215686}
\definecolor{c7}{rgb}{0.8117647058823529, 0.8117647058823529,
0.8117647058823529}

\pgfplotscreateplotcyclelist{myCycleListColor}
{%
  {blue, mark=o},
  {red, mark=square},
  {darkbrown, mark=triangle},
  {darkgreen, mark=diamond},
  {black, mark=pentagon},
  {blue, mark=*},
  {red, mark=square*},
  {darkbrown, mark=triangle*},
  {darkgreen, mark=diamond*},
  {black, mark=pentagon*},
}

\pgfplotscreateplotcyclelist{myCycleListColorNoMarkers}
{
  {black, thick},
  {blue, thick},
  {red, dotted},
  {darkbrown, loosely dashed},
  {darkgreen, loosely dotted},
}

\pgfplotsset{every axis/.append style= 
              {
                font=\small,
                mark size=2,
                line width = 0.1,
                legend style={font=\small, mark size=3, draw=none, fill=none},
                legend cell align=left,
                cycle list name=myCycleListColor,
              }
            }
            
\pgfplotstableset
{
  every odd row/.style={before row={\rowcolor[gray]{0.9}}},
  every head row/.style=
  {
    before row=
    {%
      \toprule
    },
    after row=\midrule
  },
  every last row/.style=
  {
    after row=\bottomrule
  }
}

\newif\ifdrawboundingbox

\usetikzlibrary{external} 
\tikzset{external/system call={pdflatex \tikzexternalcheckshellescape
-halt-on-error -interaction=batchmode -jobname "\image" "\texsource"}} 

\pgfkeys
{ 
  /pgf/number format/.cd,
  fixed,
  set thousands separator={\,},
  1000 sep in fractionals,
}

\newcommand{\refSection}[1]{Section~\ref{#1}}
\newcommand{\refFigure}[1]{Figure~\ref{#1}}
\newcommand{\refTable}[1]{Table~\ref{#1}}
\newcommand{\refEquation}[1]{(\ref{#1})}

\newcommand{\tensor}[1]{\boldsymbol{{#1}}}

\newcolumntype{C}[1]{>{\centering\arraybackslash}m{#1}}
\newcolumntype{R}[1]{>{\raggedright\arraybackslash}m{#1}}
\newcolumntype{L}[1]{>{\raggedleft\arraybackslash}m{#1}}

\newcommand{\delete}[1]{\xspace}

\drawboundingboxtrue
\drawboundingboxfalse

\begin{document}

\begin{frontmatter}
\title{Residual stresses in metal deposition modeling: \\ discretizations of higher order}
\author[1]{A. Özcan\corref{cor1}}
\ead{ali.oezcan@tum.de}
\cortext[cor1]{Corresponding author} 
\author[1]{S. Kollmannsberger}
\author[1]{J. Jomo}
\author[1,2]{E. Rank}
\address[1]{Chair for Computation in Engineering, Technische Universität München,  Arcisstr. 21, 80333 München, Germany}
\address[2]{Institute for Advanced Study, Technische Universität München,  Lichtenbergstraße 2a, 85748 Garching, Germany}	
\begin{abstract}
This article addresses the research question if and how the finite cell method, an embedded domain finite element method of high order, may be used in the simulation of metal deposition to harvest its computational efficiency. This application demands for the solution of a coupled thermo-elasto-plastic problem on transient meshes within which history variables need to be managed dynamically on non-boundary conforming discretizations. To this end, we propose to combine the multi-level $hp$-method and the finite cell method. The former was specifically designed to treat high-order finite element discretizations on transient meshes, while the latter offers a remedy to retain high-order convergence rates also in cases where the physical boundary does not coincide with the boundary of the discretization. We investigate the performance of the method at two analytical and one experimental benchmark. 
\end{abstract}

\end{frontmatter} 

\normalem

\tableofcontents
\newpage

\section{Introduction}

The process of metal additive manufacturing involves a highly localized heat source which moves
over a substrate. Its purpose is to change the state of the deposited metal either from powder to
liquid or from solid to liquid such that the added material bonds with the substrate. The
subsequent, rapid cooling of the heat affected zone induces undesired residual stresses, a process
well understood in welding, see for example~\cite{Goldak2005,Lindgren2007}. The numerous physical
phenomena involved in this process range from the simplified view stated above to models including
detailed weld pool dynamics~\cite{Tanaka2014}, to models resolving the micro-structure evolution in
the cooling phase (see e.g. ~\cite{Zinovieva2018,Basak2016,Korner2016}). The necessary compromise
between the complexity of the physical model and the ever increasing yet limited resources for its
numerical resolution is mainly driven by the concrete question which needs to be answered by the
simulation, (see e.g.~\cite{Lindgren2006} for a guideline in welding).

In this context, efficient numerical discretizations are desirable of which the finite element
method is the most prominent choice. A simple measure of its complexity is the number of degrees of
freedom involved in the computation. These represent the unknown coefficients of piecewise
polynomials spanned on a mesh of finite elements which resolve the highly localized gradients.
Three basic strategies are available to control the number of degrees of freedom: a)
$h$-refinement, i.e. an increase of the number of finite elements in the mesh, b) $p$-refinement,
i.e. an increase of the polynomial degree of each finite element or c) the combination of both:
$hp$-refinement. The achievable rate of convergence increases from a) to c) for problems with
locally high gradients. Unfortunately, so does the complexity of the implementation of the
respective method. This is especially the case in transient problems where the necessary adaptions
need to be kept local to the (traveling) heat affected zone.

Numerous approaches are reported in literature for locally $h$-refined discretizations, see
e.g.~\cite{Lindgren1997} as a representative of an early work on the subject in the context of
welding and~\cite{Schoinochoritis2017} for a recent review of the available variants applied to
metal additive manufacturing including a review of commercial packages for this purpose.
Refinements in $p$ beyond quadratic shape functions are scarce, but e.g.~\cite{Nubel2007}
demonstrates that it is possible to construct efficient discretizations for elasto-plastic problems
by using high order discretizations whose boundary follows the plastic front. The application of
$hp$-adaptive methods in this context is e.g. discussed in~\cite{Olesky2015} where exponential
convergence rates were achieved for boundary-conforming discretization. Recently, high-order $h$
refinements have also become of interest in the context of Isogeometric Analysis, see
e.g.~\cite{Elguedj2014}. All these publications treat boundary conforming discretizations.
Extensions to non-boundary conforming finite element discretizations are presented
in~\cite{Abedian2014,Taghipour2018} in the context of the finite cell method. Therein, it is
demonstrated that high-order convergence rates may also be achieved even if the physical boundaries
of the domain do not coincide with the discretization of the mesh. This is a desirable feature in
the application at hand where the physical domain grows with time. Additionally, highly localized
gradients need to be followed dynamically though the course of the simulation.

As a remedy for an accurate resolution of transient gradients on non-boundary conforming domains we
advocate a combination of both, the multi-level $hp$-method and the finite cell method. To this
end, we first introduce these numerical methods in~\cref{sec:numericalMethods} where we begin by
recalling the multi-level $hp$-method. It offers a relatively simple management of the degrees of
freedom for transient discretizations of $hp$-type. Additionally, we present some novel but
straight-forward extensions necessary for elasto-plastic computations on transient discretizations
of $hp$-type. We then present the computation of elasto-plasticity in the framework of the finite
cell method~\cite{Duster2008} and proceed to a short description of the thermo-elasto-plastic
problem in~\cref{sec:elastoplasticity}\,.

After the basics of the numerical methods are presented, we evaluate their performance against
three benchmark examples. In~\cref{sec:plasticSphere} a sphere under internal pressure in which a
plastic front evolves is investigated. While the exact location of the plastic front is known in this specific case, 
in a more general setting it would not be unknown. We, thus, chose a grid-like discretization whose boundaries neither coincide with the
evolving plastic front nor with the physical boundaries of the problem setup itself. We demonstrate
that it is possible to capture the plastic front and the stress states on non-boundary conforming
domains using higher order $h-$ refinements and obtain higher order convergence rates. We then move
to the thermo-elasto-plastic setting in~\cref{sec:elastoplasticBar} to make the point that it is
also possible to capture the plastic stress states in a coupled setting accurately. The final
numerical example presented in~\cref{sec:metalDepositionExperiment} is chosen to demonstrate that
the presented methodology is also capable of reproducing stress states of real experiments. We
conclude the article by pointing out the potential and limits of the presented approach in~\cref{sec:conclusions}\,.
\section{Numerical Methods} \label{sec:numericalMethods}
This section serves to introduce and further develop the numerical techniques which are used and
evaluated in~\cref{sec:numericalExamples}\,. To this end, the multi-level $hp$-method is first
presented in~\cref{sec:multilevelhp} followed by a description on the finite cell
method~\cref{sec::fcm}\,. They are then combined to solve problems in elasto-plasticity as introduced
in~\cref{sec:elastoplasticity}\,. All these sections form the background for solving
thermo-elasto-plastic problems as introduced in~\cref{sec:coupling}\,.

\subsection{The multi-level hp-method}
\label{sec:multilevelhp}
\begin{figure}[t]
  \begin{minipage}[b]{.44\textwidth}
    \begin{center}
    \includegraphics[width=0.8\textwidth]{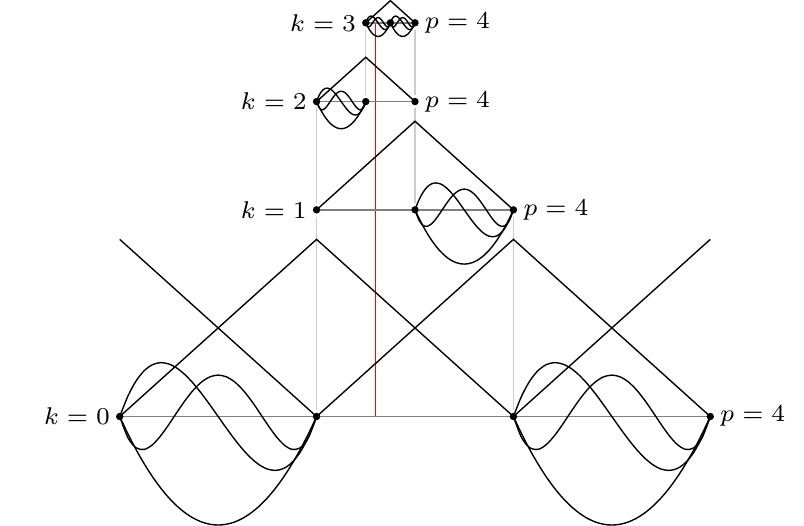}%
    \end{center}
    \begin{center}
      \small
      \vspace{-0.3cm}
      (a) One-dimensional case
    \end{center}
    \begin{center}
    \includegraphics[width=0.75\textwidth]{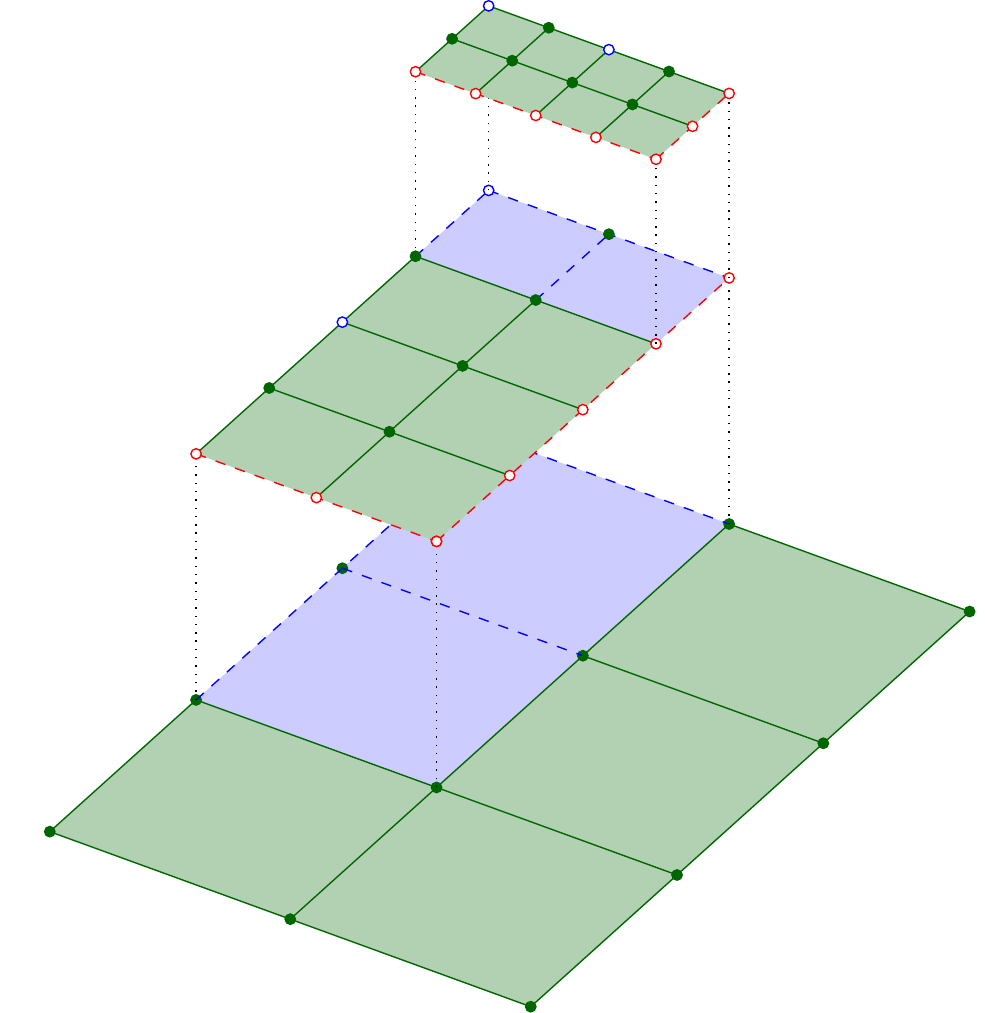}%
      \small \\
      (b) Two-dimensional case
    \end{center}
  \end{minipage}%
  \hfill%
  \begin{minipage}[b]{.55\textwidth}
    \includegraphics[width=0.88\textwidth]{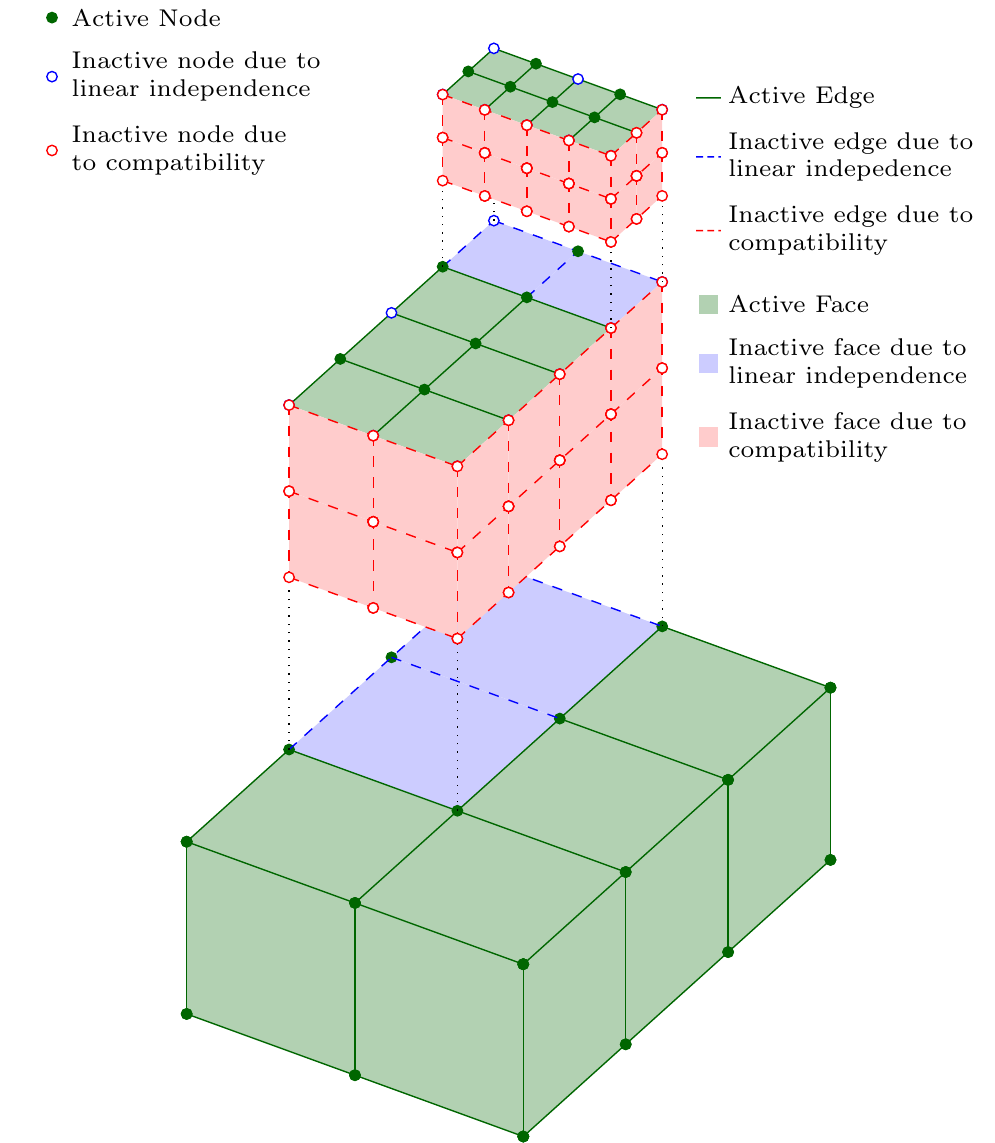}%
    \begin{center} 
        \small
        (c) Three-dimensional case
      \end{center}
  \end{minipage}
  \caption{Main idea of the multi-level hp-method in different spatial dimensions \cite{Zander2016}.\label{fig::multilevelhpIdea}}
\end{figure} 

The aforementioned multi-level hp-method is a powerful scheme for performing local mesh adaptation in an efficient manner. It was first introduced in~\cite{Zander2015}
and aims at a simple degree of freedom management for dynamic discretizations of high order by performing local refinements based on superposition. 
For this purpose, an initial discretization of coarse base elements is overlayed locally with multiple layers of finer overlay elements so as to better capture
the solution behavior such as locally high gradients, see~\cref{fig::multilevelhpIdea}\,.
This is in contrast to standard $hp$-methods which perform refinement by replacement where coarse elements are replaced by a set of smaller elements. 

The methodology to overlay finite elements as in the multi-level $hp$-method requires the enforcement of linear independence of the basis functions and  
compatibility of the ansatz space at the boundary of the refinement zone. This is achieved in a straightforward manner by leveraging the direct association of topological components 
(nodes, edges, faces, volumes) with the degrees of freedom. Compatibility and linear independence are ensured through the deactivation of specific topological
components. This deactivation is governed by simple rule-sets, which work for different spatial dimensions \cite{Zander2015,Zander2016} alike. 
This forms a key strength of the method as compared to classic $hp$-methods as the simplicity of the rule-sets allow for an easy treatment of dynamic meshes 
with arbitrary levels of hanging nodes. The multi-level $hp$-method results in $C^0$-continuity over
the complete refinement hierarchy in base elements while $C^{\infty}$-continuity is obtained in the 
leaf elements --- elements on the highest refinement level with no children.

\subsection{The finite cell method}
\label{sec::fcm}
\begin{figure}[h!]
  \begin{center}
    \def\svgwidth{\columnwidth}
    \resizebox{0.9\textwidth}{!}{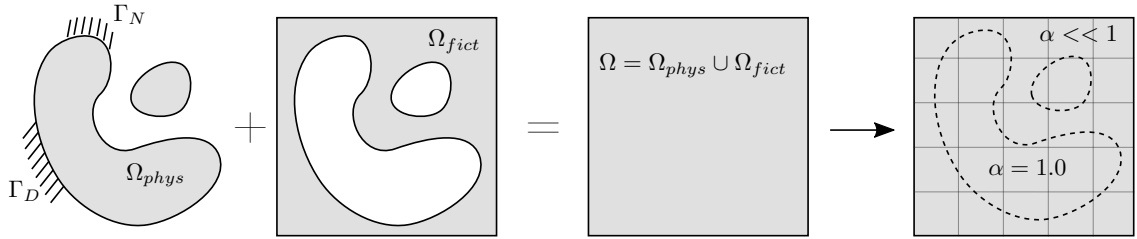}
  \caption{Core idea of the finite cell method.}
  \label{fig::finitecellmethod} 
  \end{center} 
\end{figure}

The ultimate in the scope of this article is to provide a framework which leaves as much geometric
and topological freedom for the emerging additively manufactured artifact as possible. At the same
time, the effort for mesh generation should be kept as low as possible. These are two main features
of immersed methods, a class of advanced discretization techniques that significantly reduce the
effort of mesh generation by utilizing a non-boundary conforming domain discretization. They have,
therefore, emerged as the method of choice to perform numerical simulations on bodies with a
complex shape, topology or a combination of both. The finite cell method (FCM) introduced
in~\cite{Parvizian2007,Duster2017}, is a prominent representative of immersed methods whose core idea is to
combine the advantages of the fictitious domain approach with the computational efficiency of
discretizations of high order. Its basic idea is depicted in~\cref{fig::finitecellmethod}\,. A body of complex shape and
topology defined on a physical domain $\Omega_{\text{phys}}$ is extended by a fictitious domain $\Omega_{\text{fict}}$. Their union yields a
computational domain $\Omega = \Omega_{\text{phys}} \cup \, \Omega_{\text{fict}} $ with a simple
boundary which can be discretized using well shaped finite elements. These form the support of
the ansatz functions and are termed finite cells as their boundary is not conforming to the
boundary of the original, physical boundary $\partial\Omega_{\text{phys}}$. The physical domain
must then be recovered at the level of the numerical integration of the associated bi-linear and
linear forms. An indicator function $\alpha$ is introduced to classify points belonging to the
physical or fictitious domain. Points within the physical domain are assigned a value $\alpha = 1$,
whereas points in the fictitious domain have a value $\alpha \ll 1$.
While a detailed description of the FCM can be found in~\cite{Parvizian2007,Duster2008,Duster2017}, the short and
concise rehearsal of its basic ideas above suffices to set the stage for its extension to elasto-plastic problems as laid out in the next section.

\subsection{Elasto-plasticity with the multi-level hp method} \label{sec:elastoplasticity}
The stage is now set to develop an efficient elasto-plastic formulation which combines the
multi-level $hp$-method with the finite cell method as presented in~\cref{sec:multilevelhp,sec::fcm}\,.

The finite cell method has already been successfully used in the field of plasticity. The work
most relevant to the paper at hand is published in~\cite{Abedian2013a} in which the FCM was
implemented for the $J_2$ flow theory with nonlinear, isotropic hardening for small displacements and
small strains. Further, in~\cite{Abedian2014} it is demonstrated that the FCM leads to more
efficient discretizations than the standard, boundary conforming $h$-version of the finite element
method delivers. The FCM has recently also been extended to nearly incompressible finite strain
plasticity with complex geometries, see~\cite{Taghipour2018}. All these investigations were carried
out on static discretizations in the sense that no dynamic refinement was applied locally to
capture transient plastic fronts. However, in layered deposition modeling the size of the traveling
heat source is comparatively small w.r.t. the rest of the computational domain. In these
applications, a dynamic refinement and coarsening offers a discretizational advantage because the
computational effort can thereby be balanced with the accuracy needed locally. However, this
requires a dynamic management of primal and internal variables.

To facilitate a comprehensive but concise presentation of this subject, the next sections start
with the classic formulation of $J_2$ plasticity in~\ref{sec:classicL2}\,. The standard setting is
then cast into the FCM formalism in~\cref{sec:elastoplasticfcm} before the integration of the
multi-level $hp$-method into the elasto-plastic FCM framework is treated
in~\cref{sec:transferOfVariables} along with the associated management of primal and internal variables.

\subsubsection{Classic $J_2$ plasticity} \label{sec:classicL2}
The classic weak form of equilibrium in a solid body is given by
\begin{equation} 
G(\textbf{\textit{u}}, \tensor{\eta}) = \int\limits_\Omega
\tensor{\varepsilon}(\tensor{\eta}) : \tensor{\sigma} \,\mathrm{d}\Omega -
\int\limits_\Omega
\tensor{\eta} \cdot  \textbf{\textit{b}} \,\mathrm{d}\Omega -
\int\limits_{\Gamma_N}  \tensor{\eta} \cdot 
\textbf{\textit{t}} \,\mathrm{d}\Gamma = 0,
\label{eq::weakFormOfMechanical}
\end{equation}
where $\textbf{\textit{u}}$ is the displacement field, $\tensor{\eta}$ is the
test function, $\tensor{\sigma}$ the Cauchy stress tensor, $\tensor{\varepsilon}$ the strain tensor, $\textbf{\textit{b}}$ are the body forces acting on the domain,
$\Omega$, and $\textbf{\textit{t}}$ is the prescribed traction on the Neumann
boundary, $\Gamma_N$.
When a nonlinear material is utilized, the stress state is not only a
function of the instantaneous strain. Instead, the stress state also depends on the history of
the loads the body was subjected to. Consequently, the weak form
\refEquation{eq::weakFormOfMechanical} becomes nonlinear and is solved
incrementally for each time step $[t_n,t_{n+1}]$. Within each time step the
internal variables, $\tensor{\lambda}$, which contain the history of the
material, are assumed constant. Linearization of the weak form with respect to
the unknown $\textbf{\textit{u}}$ around $\textbf{\textit{u}}_{n+1}^{(i)}$ , which is
the solution at iteration i, is given e.g. in~\cite{de2009computational} and reads 
\begin{equation} 
\begin{aligned} 
\int\limits_\Omega
\tensor{\varepsilon}(\tensor{\eta}) : \tensor{D} :
\tensor{\varepsilon}(\delta \textbf{\textit{u}}) \,\mathrm{d}\Omega =
- \int\limits_\Omega
\tensor{\varepsilon}(\tensor{\eta}) :
\tensor{\sigma}_{n+1}(\tensor{\lambda}_n,\tensor{\varepsilon}(\textbf{\textit{u}}_{n+1}^{(i)}))
\,\mathrm{d}\Omega \\ + \int\limits_\Omega \tensor{\eta} \cdot 
\textbf{\textit{b}}_{n+1} \,\mathrm{d}\Omega + \int\limits_{\Gamma_N}  \tensor{\eta} \cdot
\textbf{\textit{t}}_{n+1} \,\mathrm{d}\Gamma.
\end{aligned}
\label{eq::linearizedWeakForm}
\end{equation}
where the fourth order tensor $\tensor{D}$ is the tangent modulus defined
as 
\begin{equation} 
\tensor{D} = \left. \frac{\partial \tensor{\sigma}_{n+1}}{\partial
\tensor{\varepsilon}}\right|_{\tensor{\varepsilon}(\textbf{\textit{u}}_{n+1}^{(i)})}.
\label{eq::tangentModulusTensor}
\end{equation}

\subsubsection{The elasto-plastic finite cell method} \label{sec:elastoplasticfcm}
In order to apply the finite cell method, the domain integrals in \cref{eq::linearizedWeakForm} are multiplied with the indicator function
$\alpha$ such that
\begin{equation} 
\begin{aligned}
\int\limits_\Omega
\alpha\,\tensor{\varepsilon}(\tensor{\eta}) : \tensor{D} :
\tensor{\varepsilon}(\delta \textbf{\textit{u}}) \,\mathrm{d}\Omega =
- \int\limits_\Omega \alpha\,\tensor{\varepsilon}(\tensor{\eta}) :
\tensor{\sigma}_{n+1}(\tensor{\lambda}_n,\tensor{\varepsilon}(\textbf{\textit{u}}_{n+1}^{(i)}))
\,\mathrm{d}\Omega \\ + \int\limits_\Omega \alpha\,\tensor{\eta} \cdot 
\textbf{\textit{b}}_{n+1} \,\mathrm{d}\Omega + 
\int\limits_{\Gamma_N}  \tensor{\eta} \cdot \bar{\textbf{\textit{t}}}_{n+1}
\,\mathrm{d}\Gamma.
\end{aligned}
\label{eq::linearizedWeakFormAlpha}
\end{equation}
Since the solution in the fictitious domain is unphysical, computing the
tangent modulus or stresses in the fictitious domain causes unnecessary computational overhead.
Therefore, the deformation in the fictitious domain is neglected and stresses are assumed
to be zero ($\tensor{\sigma}_{fict} = \tensor{0}$). Moreover, the tangent
modulus is taken as the elastic tangent ($\tensor{D}_{fict} = \tensor{D}^e$).
Incorporating these assumptions into~\cref{eq::linearizedWeakFormAlpha} provides the following
linearized weak form of equilibrium for elasto-plastic problems
\begin{equation} 
\begin{aligned} 
\int\limits_{\Omega_\text{phys}}
\tensor{\varepsilon}(\tensor{\eta}) : \tensor{D} :
\tensor{\varepsilon}(\delta \textbf{\textit{u}}) \,\mathrm{d}\Omega +
\int\limits_{\Omega_\text{fict}}
\alpha\,\tensor{\varepsilon}(\tensor{\eta}) : \tensor{D}^{e} :
\tensor{\varepsilon}(\delta \textbf{\textit{u}}) \,\mathrm{d}\Omega = \\ 
- \int\limits_{\Omega_\text{phys}}
\tensor{\varepsilon}(\tensor{\eta}) :
\tensor{\sigma}_{n+1}(\tensor{\lambda}_n,\tensor{\varepsilon}(\textbf{\textit{u}}_{n+1}^{(i)}))
\,\mathrm{d}\Omega + 
\int\limits_\Omega \alpha\,\tensor{\eta} \cdot  \textbf{\textit{b}}_{n+1}
\,\mathrm{d}\Omega + 
\int\limits_{\Gamma_N}  \tensor{\eta} \cdot \bar{\textbf{\textit{t}}}_{n+1}
\,\mathrm{d}\Gamma.
\end{aligned}
\label{eq::linearizedWeakFormFCM}
\end{equation}

\begin{figure}
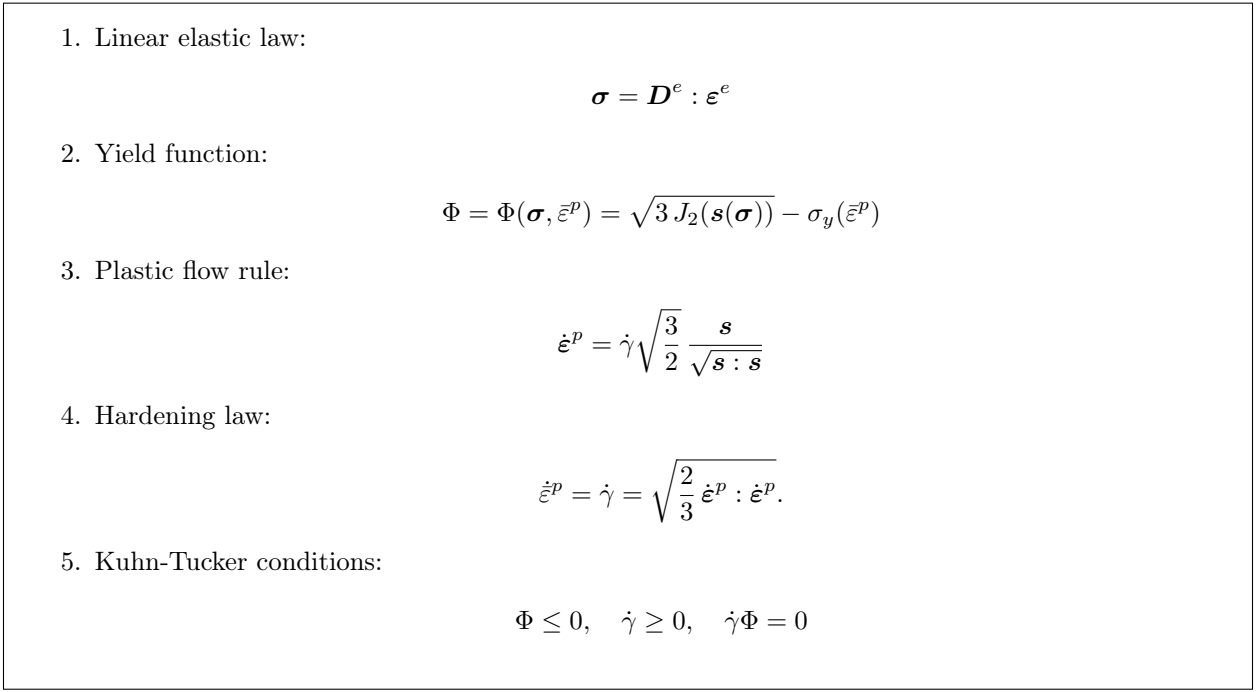

\begin{center}
\begin{framed}
\begin{enumerate}
\item Linear elastic law: 
      \begin{align*}
        \tensor{\sigma} = \tensor{D}^e : \tensor{\varepsilon}^e
      \end{align*}   
\item Yield function:
      \begin{align*}
        \Phi = \Phi(\tensor{\sigma}, \bar{\varepsilon}^p) =
        \sqrt{3\,J_2(\tensor{s}(\tensor{\sigma}))} - \sigma_{y}(\bar{\varepsilon}^p )
      \end{align*}
\item Plastic flow rule:
      \begin{align*}
		\dot{\tensor{\varepsilon}}^p =
		\dot{\gamma}\sqrt{\frac{3}{2}}\,\frac{\tensor{s}}{\sqrt{\tensor{s}:\tensor{s}}}
      \end{align*}
\item Hardening law:
      \begin{align*}
		\dot{\bar{\varepsilon}}^p = \dot{\gamma} = \sqrt{\frac{2}{3}\,
		\dot{\tensor{\varepsilon}}^p : \dot{\tensor{\varepsilon}}^p}.
      \end{align*}            
\item Kuhn-Tucker conditions:
      \begin{align*}
        \Phi \leq 0, \quad \dot{\gamma} \geq 0, \quad \dot{\gamma}\Phi = 0
      \end{align*}
\end{enumerate}
\end{framed}
\caption{Summary of the rate-independent von Mises associative model with nonlinear isotropic
hardening}
\label{fig::vonMisesModel}
\end{center}
\end{figure}

In this paper, a rate-independent von Mises plasticity model with
nonlinear isotropic hardening is considered whereby the displacements and
strains are assumed to be small. The classic, additive decomposition of the
strain tensor into an elastic and a plastic counterpart is then applied such that
\begin{equation} 
\tensor{\varepsilon} = \tensor{\varepsilon}^e + \tensor{\varepsilon}^p.
\label{eq::additiveDecomposition3D}
\end{equation}
\refFigure{fig::vonMisesModel} summarizes this model, in which the
internal variables $\tensor{\lambda}$ are given by
\begin{equation} 
\tensor{\lambda} = \{ \,\tensor{\varepsilon}^p, \bar{\varepsilon}^p \, \},
\label{eq::internalVariablesSetVonMises}
\end{equation}
where $\bar{\varepsilon}^p$ is the equivalent plastic strain. 

\subsubsection{Multi-level $hp$-adaptivity and the transfer of primary and internal variables} \label{sec:transferOfVariables}
In the case of dynamic multi-level \textit{hp}-refinements the
primary variables given by the displacement field, $\textbf{\textit{u}}$, and internal variables
given by $\tensor{\lambda}$, need to be transferred from the old discretization prior to the
refinement to the new discretization after the refinement was carried out. Since the displacement
field over the domain is discretized by the basis functions, a global $C^0$ continuous description
is available which is directly transferred to the new discretization by means of a global
$L_2$-projection. However, in the general elasto-plastic finite element procedure, the evolution of
the internal variables is computed at \emph{local} integration points via the plastic flow rule and
the hardening law. As no $C^0$ continuous discretization is readily available, the transfer of
these variables from the integration points of the old discretization to the integration points of
the new one is more involved.

As a remedy, several strategies have been developed in literature. The simplest is the point-wise
transfer of the internal variables~\cite{Peric1996}. Therein, an area is associated to each old
integration point within which constant values of the history variables are assumed. This leads to
a discontinuous approximation of the history variables within the concerned finite elements.
Another strategy is the element-wise transfer \cite{Ortiz1991}, where the internal variables are
interpolated by local functions. This strategy can result in an approximation with discontinuities
over element boundaries. A cure is offered by the large group of strategies involving nodal
projections. Therein the internal variables are first transferred from integration points of the
old discretization to the nodal degrees of freedom of that old discretization. This is achieved
in~\cite{Peric1996} and \cite{Lee1994} by extrapolating the values from integration points to the
nodal dofs for each element locally and later averaging these extrapolated values.
Another possibility is to use the super convergent patch recovery technique as presented in
\cite{Zienkiewicz1992a, Zienkiewicz1992b,Boroomand1999,Khoei2009}, for the transfer of the values
from the integration points of the old mesh to its nodes. In this approach a patch consists of
elements that surround a node and the internal variables are fitted to continuous polynomials over
each patch by a least squares method and are then interpolated to the nodal points. Independent of
how the values are obtained at the nodes of the old discretization, the next step is to interpolate
them to the new nodal points in the new mesh. Finally, the values of the internal variables at the
new integration points are interpolated from the values of the new nodal points. Since this
strategy acts as a smoothing operator, the localized internal variables are smeared out over a
larger part of the domain. This can be a drawback in elasto-plastic analysis, because plasticity is
a local effect in most cases.

Unfortunately, none of these techniques are suitable for the multi-level \textit{hp}-adaptivity.
This is mainly due to the fact that the high order basis functions are modal functions, which are
not associated to specific nodes. For this reason a modified version of the element-wise transfer
\cite{Ortiz1991} is advocated in the paper at hand in which each component of the internal
variables, $\lambda_i$, is approximated such that
\begin{equation} 
\lambda_i \approx \lambda^h_i = \mathbf{P}(r,s,t)\,\mathbf{c}
\label{eq::approximationOfInternalVariables}
\end{equation}
where $\mathbf{P}$ is a vector of integrated Legendre polynomials as used in
\textit{p}-version of FEM \cite{Duster2001} and $\mathbf{c}$ is the
corresponding coefficient vector. It is explicitly pointed out that these polynomials
only span the leaf element, where $(r,s,t)$ denotes the local coordinates of
this element. The coefficient vector, $\mathbf{c}$ is determined by applying
a least square fit to the discrete integration point values,
$\lambda^d_i(r_j,s_j,t_j)$, inside the leaf element. This is done by minimizing the function, $F(\mathbf{c})$,
\begin{equation} 
F(\mathbf{c}) = \sum_{j=1}^{n_{gp}}
\big(\lambda^d_i(r_j,s_j,t_j)-\mathbf{P}(r_j,s_j,t_j)\,\mathbf{c}\big)^2,
\label{eq::leastSquaresFunction}
\end{equation}
which is carried out by differentiating $F$ with respect to $\mathbf{c}$ as
\begin{equation} 
\frac{\partial F}{\partial \mathbf{c}} = \mathbf{b} - \mathbf{A}\, \mathbf{c} =
0,
\label{eq::leastSquaresMinimization}
\end{equation}
with
\begin{equation} 
\mathbf{A} =
\sum_{j=1}^{n_{gp}}\mathbf{P}^T(r_j,s_j,t_j)\,\mathbf{P}(r_j,s_j,t_j), \quad
\mathbf{b} = \sum_{j=1}^{n_{gp}}
\mathbf{P}^T(r_j,s_j,t_j) \, \lambda^d_i(r_j,s_j,t_j).
\end{equation}
Equation \refEquation{eq::leastSquaresMinimization} leads to the following
linear system of equations 
\begin{equation} 
\mathbf{A}\, \mathbf{c} = \mathbf{b},
\end{equation}
which needs to be solved for each component of the internal variables.

\begin{figure}[htb]
  \begin{center}
       \includegraphics[]{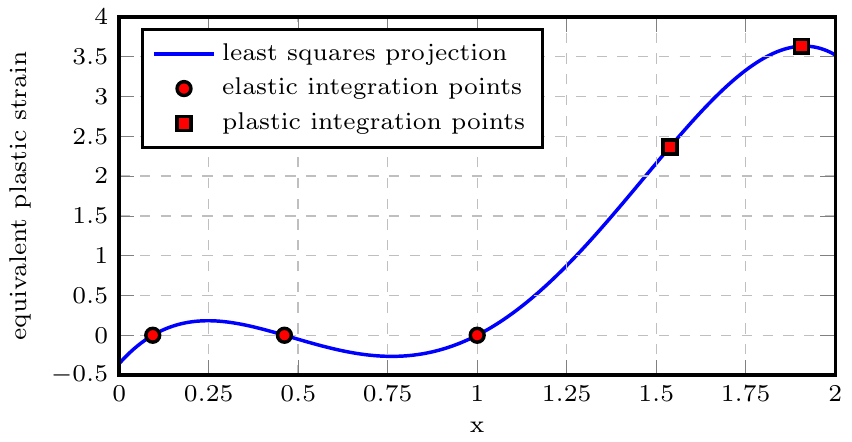}
  \end{center}
  \caption{Convergence}
  \label{fig::oscilatingProjection}
\end{figure}

If this projection strategy is applied to the elements that
contain the plastic front, i.e. in which only a group of integration points in
the element accumulated plastic strains, then the proposed least squares approximation may lead to
oscillations and spurious values. In some cases it is even possible to obtain invalid negative values for the equivalent
plastic strain. \refFigure{fig::oscilatingProjection} illustrates this for the polynomial order $p=4$ for a simple one dimensional example in which the equivalent plastic strain should be nonzero for the two plastified integration points and zero for the three elastic integration points. This function
cannot be represented by polynomials and the well known Gibbs phenomenon occurs.

\begin{figure}[htb]%
  \begin{center}%
    \subfloat[Interpolation cases]%
    {     
      \includegraphics[]{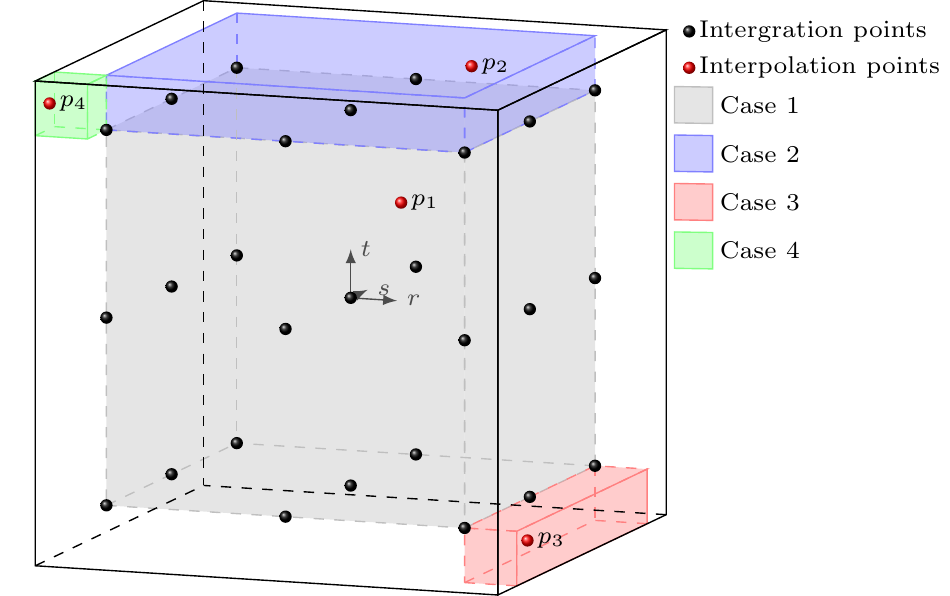}
      \label{fig::interpolationZones}
    }
    \hfill%
    \subfloat[Case 1: Trilinear interpolation]%
    {     
      \includegraphics[]{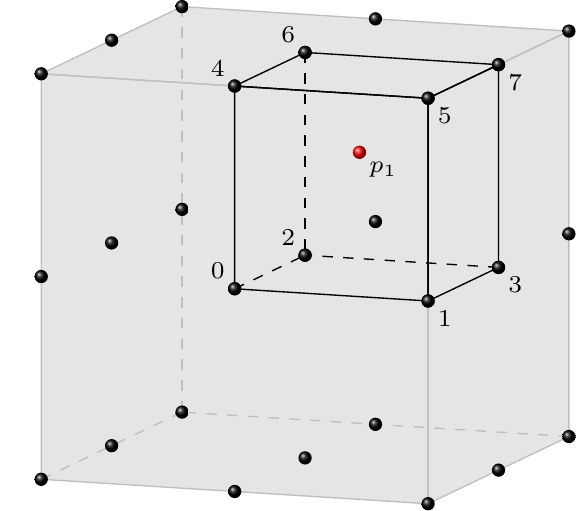}%
      \label{fig::trilinearInterpolation}
    }
    \hfill%
    \subfloat[Case 2: Bilinear interpolation]%
    {
      \includegraphics[]{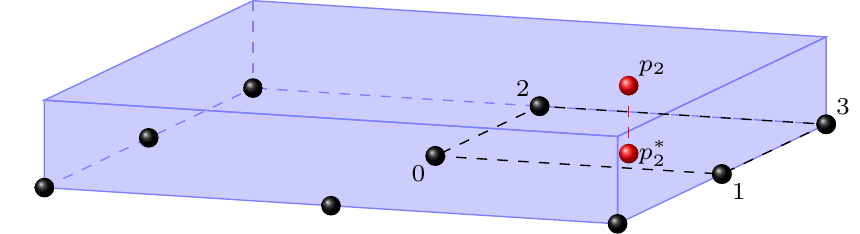}%
      \label{fig::bilinearInterpolation}
    }
    \hfill%
    \subfloat[Case 3: Linear interpolation]%
    {
     \includegraphics[]{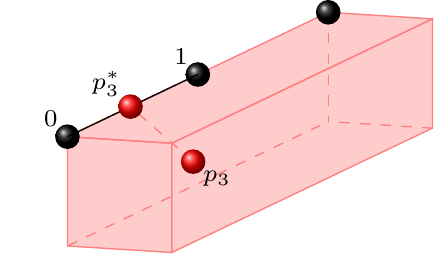}%
     \label{fig::linearInterpolation} 
    }
    \caption{Interpolation in interface element}%
    \label{fig::interpolationInInterfaceElement}%
  \end{center}%
\end{figure}

In order to circumvent this problem, the values in those elements containing a elastic-plastic interface are
not projected with the least squares method. Instead, for these elements we advocate the following approach depicted in~\refFigure{fig::interpolationZones}\,. 
Four cases are distinguished and characterized by the location of the integration points w.r.t. the
concerned interpolation point. The first case is given in \refFigure{fig::trilinearInterpolation}
where the interpolation point $p_1$ is located inside the gray bounding box defined by all
integration points. In this case the internal variables, $\lambda_i$ at $p_1$ are approximated by a
trilinear interpolation from the values of the eight integration points closely surrounding it. In
all the other cases the interpolation point $p$ is located outside the gray bounding box.
Therefore, its projection onto the corresponding surface, edge or corner of the bounding box,
$p^*$, is used for interpolation. \refFigure{fig::bilinearInterpolation} depicts the second case in
which the projection of $p_2$ is surrounded by only four points. These are then used for bilinear
interpolation. The third case occurs when the projected point is surrounded by only two points as
demonstrated in \refFigure{fig::linearInterpolation}\,. In this case linear interpolation is applied.
The last case occurs when the original interpolation point lies in the corner quadrants of the
finite cell. In this case values at the closest integration point are used. We will demonstrate
in~\cref{sec:numericalExamples} that this approach is feasible for dynamic multi-level
$hp$-discretizations in an elasto-plastic setting.

\subsection{Coupling to thermal problems} \label{sec:coupling}
Metal deposition is a multi-physics problem,
where the highest temperature gradients as well as the phase changes occur in
close vicinity of the moving laser beam. The phase changes between liquid and
solid states are simulated by the model introduced by Celentano et.
al. \cite{Celentano1994} which uses the discretized weak form
of equation \refEquation{eq::thermal}. Further details regarding this
method and its application to thermal analysis of the selective laser melting
process in the framework of multi-level \textit{hp}-adaptivity can be found in
~\cite{Kollmannsberger2016}.
\begin{equation}
\rho c \frac{\partial T}{\partial t} +  \rho L \frac{\partial f_{pc}}{\partial t} -
\nabla \cdot (k\nabla T) = Q
\label{eq::thermal}
\end{equation}

A dynamically adaptable data container (multi-level grid) is used
to keep track of the physical domain $\Omega_{phys}$ during the deposition
process. The container represents a dynamic octree. It stores the current
material state of a voxel at all points in time. Material states are stored in this container and not associated to finite elements. This decoupling of material and discretization facilitates the use of the finite
cell method.

Due to the changes in the temperature during the process, some regions in
the domain expand, while others contract. This generates residual stresses
in the part when the thermal load is removed and the material allowed to cool down to its
initial temperature. These stresses are computed by the quasi-static mechanical model, given
in \refSection{sec:elastoplasticity} with the addition of the thermal strain
defined as
\begin{equation}
\tensor{\varepsilon}^{th} = \gamma \,\Delta T\, \tensor{I},
\label{eq::thermal}
\end{equation}
where $\gamma$ is the thermal expansion coefficient and $\tensor{I}$ is the
second order identity tensor. Adding the thermal strain component extends the
additive decomposition of strains given in equation
\refEquation{eq::additiveDecomposition3D} to
\begin{equation} 
\tensor{\varepsilon} = \tensor{\varepsilon}^e + \tensor{\varepsilon}^p +
\tensor{\varepsilon}^{th}.
\label{eq::additiveDecompositionWithThermal}
\end{equation}
It should be noted that, the von Mises plasticity depends solely on the deviatoric part
of the strain tensor and the thermal strain is purely hydrostatic. Therefore,
the basic structure of the elasto-plastic model introduced
in~\refSection{sec:elastoplasticity} remains as is. Only the additive
decomposition of the strains changes along with the fact that the yield stress
and the thermal expansion coefficient are now temperature dependent functions.

\begin{figure}[!ht]
\centering
    \includegraphics[]{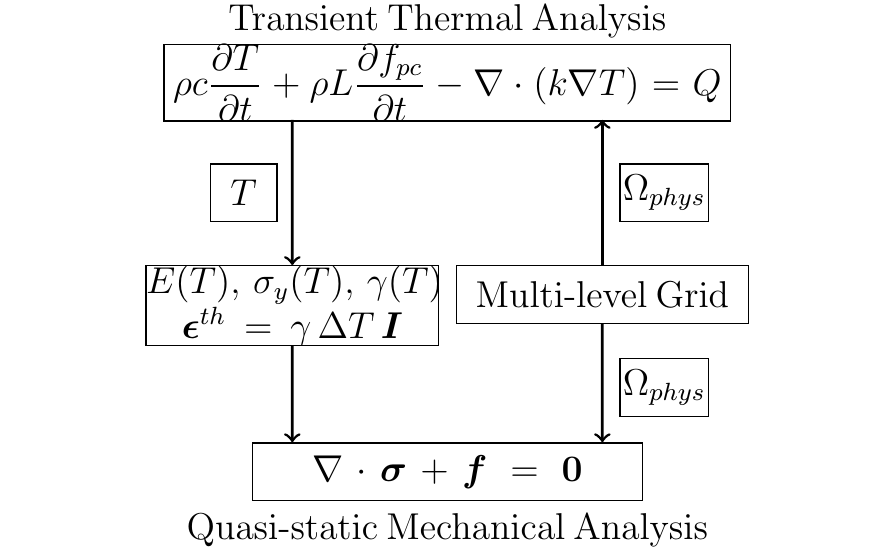}
    \caption{Thermomechanical coupling}%
    \label{fig::thermomechanicalModel}%
\end{figure}

It is assumed that displacements are small and do not produce heat, so that only
a  one-directional coupling has to be taken into account, i.e. only the
displacement field is affected by the changes in the temperature field.
As shown in \refFigure{fig::thermomechanicalModel}\,, a staggered approach is
taken for the solution of the thermomechanically coupled problem. For each time
step, the multi-level grid is first updated according to the metal deposition
such that the physical domain, $\Omega_{phys}$, is defined for both the thermal
and the mechanical problems. Then, the thermal problem is solved to obtain the
temperature distribution. Finally, the resulting temperature field is used to
compute the thermal strains and temperature-dependent material properties used to
solve the mechanical problem before the next time step is computed.

\section{Numerical Examples} \label{sec:numericalExamples}
All presented examples have benchmark character and are chosen thoroughly to test the main aspects
of the proposed methodology. The first example discussed in~\cref{sec:plasticSphere} is chosen to
test if higher order convergence rates are possible in elasto-plastic computations and under what
circumstances they decay. To this end, a dynamic multi-level $hp$-refinement as presented
in~\cref{sec:multilevelhp} is carried out towards the plastic front in several load steps whereby
the plastic front itself is not resolved in a boundary conforming manner but travels through finite
cells. This situation is expected for engineering applications. The example also tests the finite
cell formulation presented in~\cref{sec::fcm} in the elasto-plastic embedded domain setting as
described in~\cref{sec:elastoplasticity} by resolving the boundary of the physical domain only on
integration level. It also serves as a test for the transfer of history variables as detailed in~\cref{sec:transferOfVariables}\,.

The second example given in~\cref{sec:elastoplasticBar} is chosen to verify the
thermo-elasto-plastic coupling procedure laid out in~\cref{sec:coupling}\,. Both fields are
discretized on separate computational grids which are refined and de-refined individually according
to the requirements of the corresponding field variable.

The final example presented in~\cref{sec:metalDepositionExperiment} tests the combination of the
methodology on a experimental benchmark used in welding. As such, it constitutes a first step
towards a verification of the proposed methodology for layer deposition modeling.

\subsection{Internally pressurized spherical shell} \label{sec:plasticSphere}
\begin{figure}[ht]%
  \begin{center}%
    \subfloat[An octant of the spherical shell]%
    {
      \includegraphics[width=0.48\textwidth]{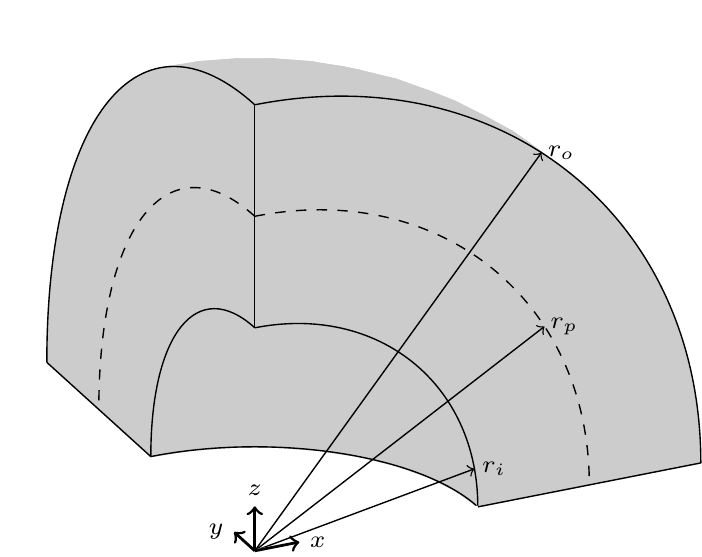}
      \label{fig::internallyPressurizesSphereSetup}
    }
    \hfill%
    \subfloat[Smart octree integration on $4\times4\times4$ base mesh]%
    {
      \includegraphics[width=0.45\textwidth]{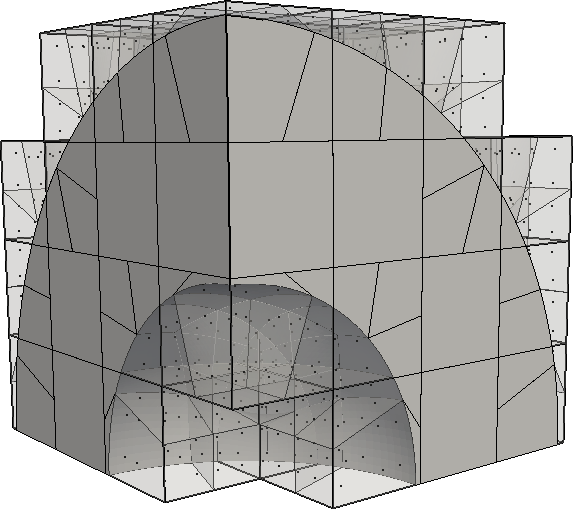}
      \label{fig::sphereSmartOctree}
    }
    \caption{Set up of internally pressurized sphere}%
    \label{fig::internallyPressurizesSphere}%
  \end{center}%
\end{figure}

In this example internal pressure, $P$, is applied in increments to a spherical
shell, which is composed of an elastic perfectly plastic material with a von
Mises yield criterion. \refFigure{fig::internallyPressurizesSphereSetup}
illustrates the problem setup, where $r_i$, $r_o$ and $r_p$ are inner radius,
outer radius and the radius of the plastic front, respectively.
Hill \cite{Hill1998} provides an analytical solution to this problem. Equation
\refEquation{eq::plasticFrontEquation} gives the relation between the location
of the plastic front, $r_p$, and the applied internal pressure, $P$, wherein
$\sigma_y$ is the yield stress.
\begin{equation}
P = 2\,\sigma_y \ln\bigg( \dfrac{r_p}{r_i}\bigg) + \dfrac{2\,\sigma_y}{3}\bigg(
1- \dfrac{r_p^3}{r_o^3}\bigg)
\label{eq::plasticFrontEquation}
\end{equation}
The analytical solutions of the normal stresses in spherical
coordinates are provided in equations \refEquation{eq::radialStress} and
\refEquation{eq::hoopStress}. The shear components in spherical coordinates are zero.

\begin{equation}
\sigma_{rr} =
  \begin{cases} 
   -2\,\sigma_y\bigg[ \ln\bigg( \dfrac{r_p}{r}\bigg) + \dfrac{1}{3} \bigg( 1 -
   \dfrac{r_p^3}{r_o^3}\bigg)\bigg] & \text{if } \,\,r \leq r_p \\ \\
    -\dfrac{2\,\sigma_y\,r_p^3}{3\,r_o^3} \bigg( \dfrac{r_o^3}{r^3} - 1  \bigg) 
    & \text{if } \,\, r > r_p
  \end{cases}
\label{eq::radialStress}
\end{equation}
\begin{equation}
\sigma_{\theta\theta} = \sigma_{\phi\phi} = 
  \begin{cases} 
   2\,\sigma_y\bigg[ \dfrac{1}{2} - \ln\bigg( \dfrac{r_p}{r}\bigg) -
   \dfrac{1}{3} \bigg( 1 - \dfrac{r_p^3}{r_o^3}\bigg)\bigg] & \text{if } \,\,r \leq r_p \\ \\
    \dfrac{2\,\sigma_y\,r_p^3}{3\,r_o^3} \bigg( \dfrac{r_o^3}{2\,r^3} + 1 
    \bigg) & \text{if } \,\, r > r_p 
  \end{cases}
\label{eq::hoopStress}
\end{equation}

Due to symmetry, only an octant of the whole spherical shell is considered in
the numerical model with the appropriate symmetry boundary conditions.
The shell is embedded in a Cartesian grid of size $4\times4\times4$ elements,
where finite cells located completely outside the spherical shell are removed
from the computation. For the numerical simulations the spherical shell is selected to
have an inner radius of $50$ mm and outer radius of $100$ mm. The
Young's modulus, $E$, and Poisson's ratio, $\nu$, are
set to $10$ GPa and $0.3$, respectively, and the yield stress, $\sigma_y$, is set to
$41.79389833783693$ MPa.

The geometry of the spherical shell is captured at the integration level by
utilizing the smart octree algorithm as explained in \cite{Kudela2016a}.
\refFigure{fig::sphereSmartOctree} illustrates the resulting integration points
and cells for base elements. 
\begin{figure}[htb]%
  \begin{center}%
    \subfloat[Radial stresses at P = 40 MPa]%
    {
      \includegraphics[width=0.31\textwidth]{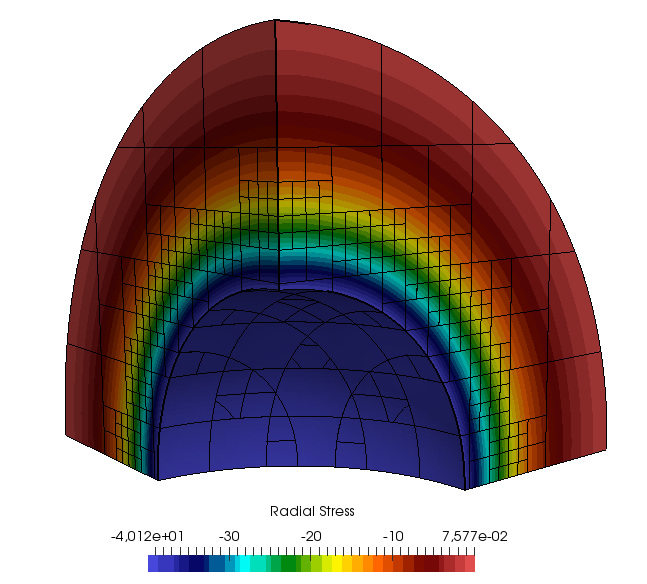}%
      \label{fig::sphereRadial40}
    }
    \hfill%
    \subfloat[Radial stresses at P = 45 MPa]%
    {
      \includegraphics[width=0.31\textwidth]{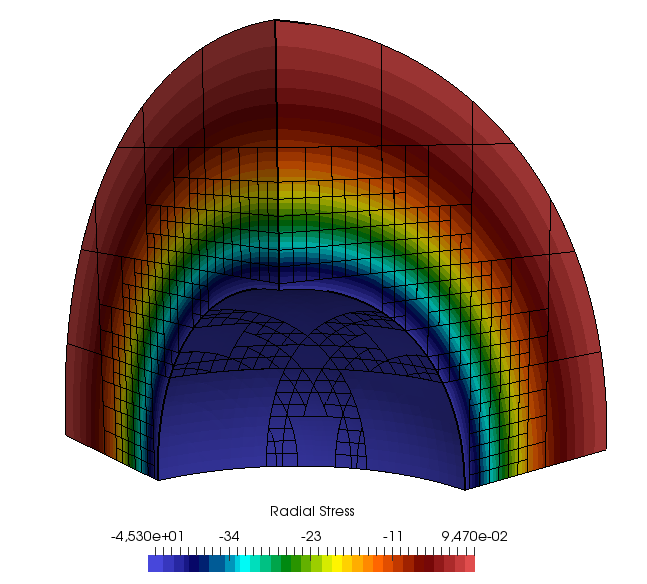}%
      \label{fig::sphereRadial45}
    }
    \hfill%
    \subfloat[Radial stresses at P = 50 MPa]%
    {
      \includegraphics[width=0.31\textwidth]{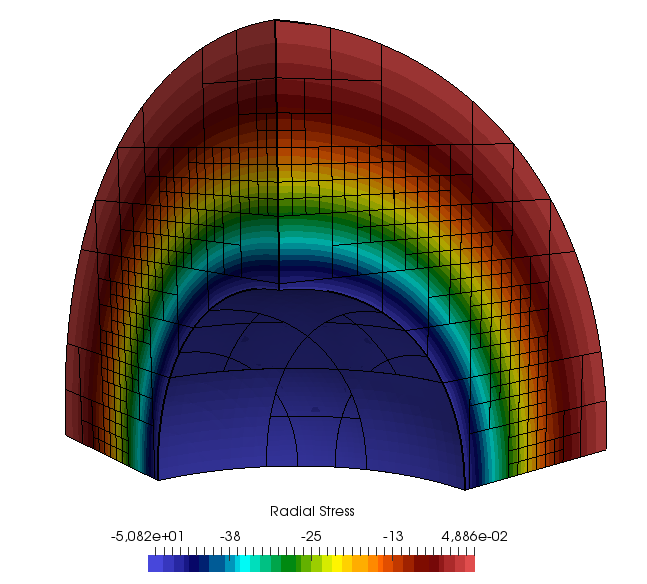}%
      \label{fig::sphereRadial50}
    }
    \caption{Radial stresses on a dynamically adaptive multi-level
    $hp$-discretization with polynomial degree $4$.}%
    \label{fig::sphereRadialResultsOnMesh}%
  \end{center}%
\end{figure}%

\refFigure{fig::sphereRadialResultsOnMesh} \,shows the distribution of the radial
stresses on a dynamically adaptive multi-level $hp$-discretization for
three load steps, where the polynomial degree of the basis functions is chosen
to be $p=4$ and maximum refinement level is set to $3$. In each step the elements
that are close to the plastic front are refined and the elements that are
further away from the plastic front are coarsened. Moreover, in each step the
primary variables related to degrees of freedom and internal variables at integration
points are transferred as explained in \refSection{sec:transferOfVariables}\,. 

\begin{figure}[ht]%
  \begin{center}%
    \subfloat[P = 40 MPa]%
    {
      \includegraphics[]{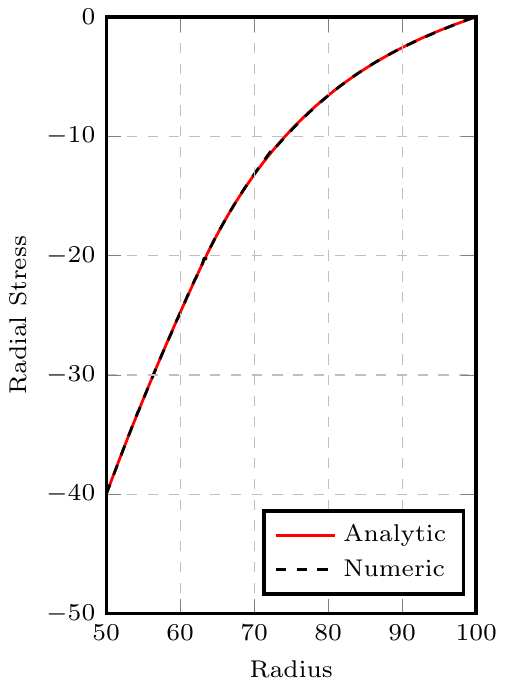}
    }
    \hfill%
    \subfloat[P = 45 MPa]%
    {
      \includegraphics[]{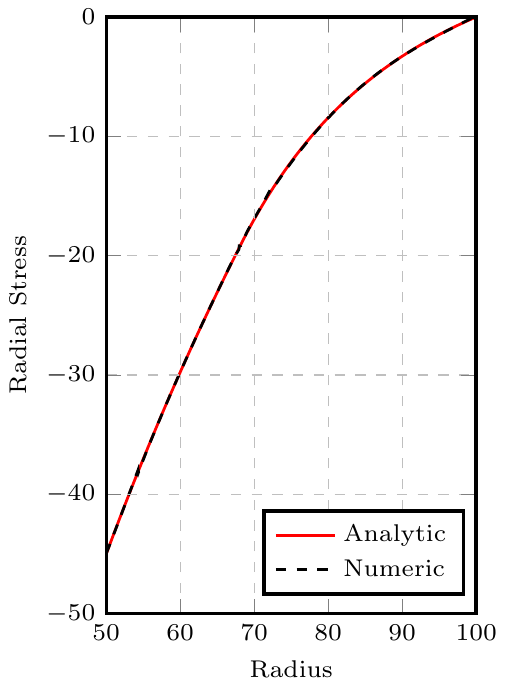}
    }
    \hfill%
    \subfloat[P = 50 MPa]%
    {
      \includegraphics[]{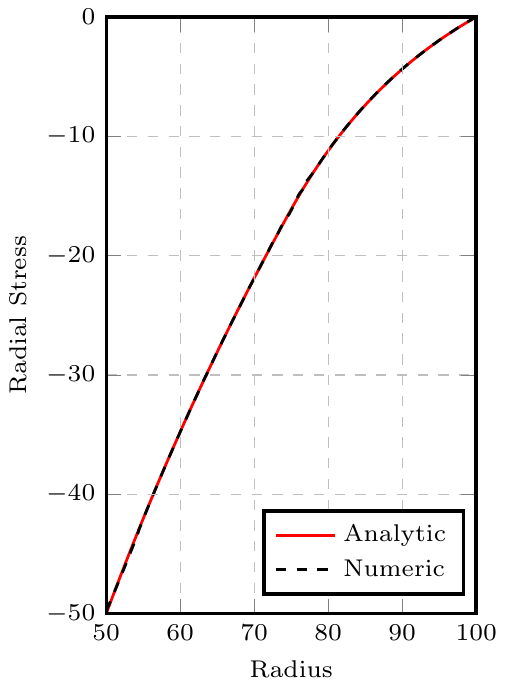}
    }
    \caption{Radial stresses along radius of the spherical shell}%
    \label{fig::radialStressesOverLoads}%
  \end{center}%
\end{figure}%

\begin{figure}[ht]%
  \begin{center}%
    \subfloat[P = 40 MPa]%
    {
      \includegraphics[]{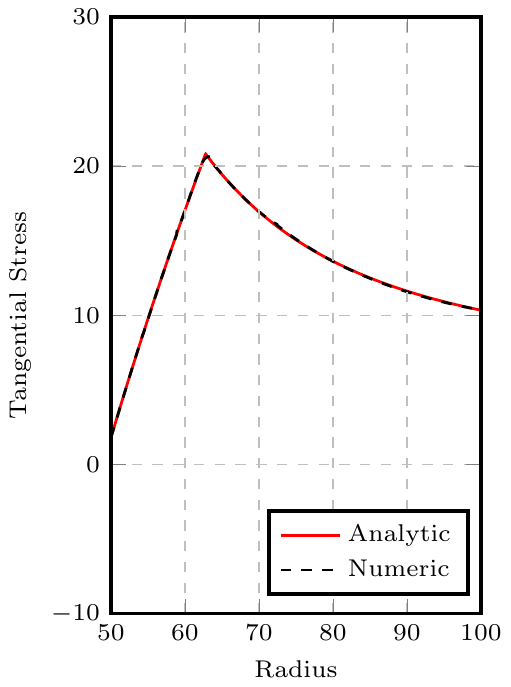}
    }
    \hfill%
    \subfloat[P = 45 MPa]%
    {
      \includegraphics[]{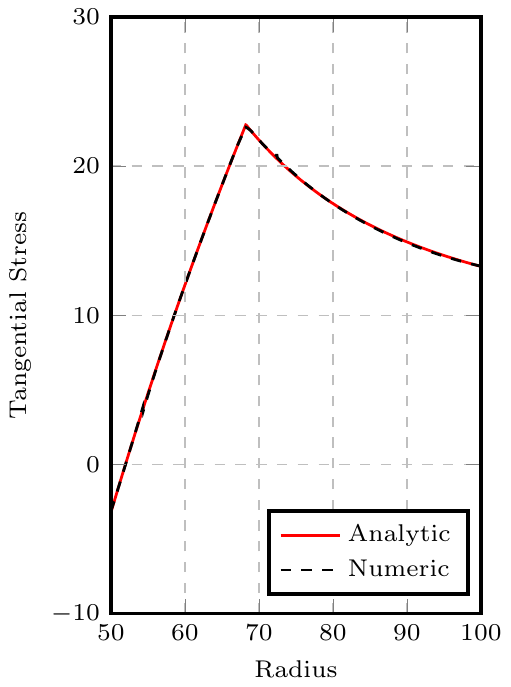}
    }
    \hfill%
    \subfloat[P = 50 MPa]%
    {
      \includegraphics[]{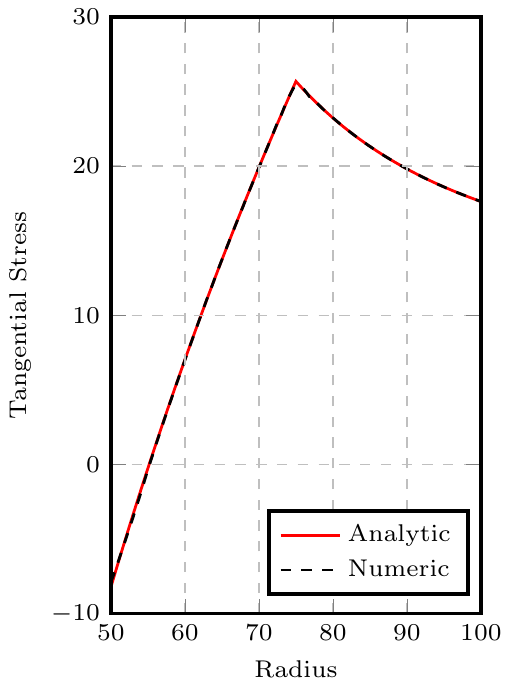}
    }
    \caption{Tangential stresses along the radius of the spherical shell}%
    \label{fig::tangentialStressesOverLoads}%
  \end{center}%
\end{figure}%

It can be seen from the equations \refEquation{eq::radialStress} and
\refEquation{eq::hoopStress} that the stresses vary only in $r$-direction and
are constant in $\theta$- and $\phi$-directions. Therefore, it is sufficient to compare the computed
stresses along the radius of the spherical shell to the analytical stresses as depicted in
\refFigure{fig::radialStressesOverLoads} and
\refFigure{fig::tangentialStressesOverLoads}\,. Both radial and tangential stress
results obtained from the numerical simulation match their analytical
counterparts for each load step. The different material behavior in elastic and
plastic regions introduces a kink in the stress field, which is more blatant for
the tangential stress. The introduced multi-level $hp$-refinement towards that elastic-plastic interface
allows the numerical solution to closely capture this behavior.

\begin{figure}[htb]
  \begin{center}
      \includegraphics[]{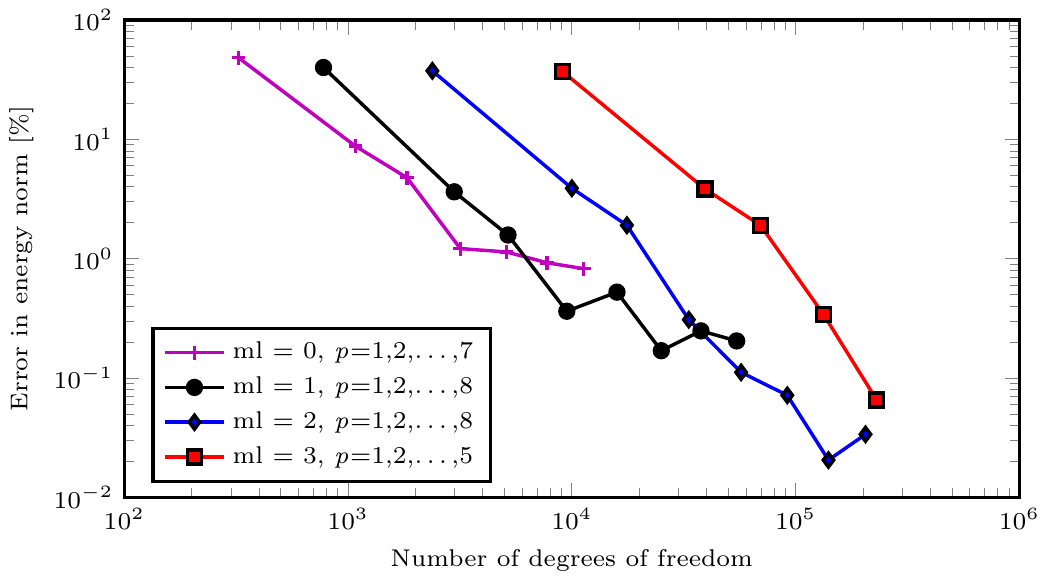}
  \end{center}
  \caption{Convergence in energy norm starting from a $4\times4\times4$ base
  mesh with multi-level $hp$ refinement (ml).}
  \label{fig::convergenceOfSphericalShell}
\end{figure}

In order to investigate the convergence properties of this elasto-plastic problem
in the framework of the multi-level $hp$-adaptive finite cell method as described
in~\cref{sec:multilevelhp,sec::fcm,sec:elastoplasticity}\,, \textit{p}-convergence studies are
performed by uniformly elevating the order \textit{p} of the polynomial shape functions, while
keeping the size of the base elements fixed. For these studies a $4\times4\times4$ base mesh is used, which
is recursively refined towards the elastic-plastic interface. Moreover, an
internal pressure $P=50$ MPa is applied in one load step for each computation.
To this end, the relative error in energy is monitored, which is defined as
\begin{equation}
e = \sqrt{ \dfrac{|U_{ex} - U_{num}|}{U_{ex}}} \times 100,
\label{eq::relativeError}
\end{equation}
with the exact internal energy, $U_{ex}$, and the numerical internal
energy, $U_{num}$, computed as 
\begin{equation}
U = \dfrac{1}{2}\int\limits_{\Omega} \tensor{\varepsilon} :
\tensor{\sigma}\,\mathrm{d}\Omega.
\label{eq::internalEnergy}
\end{equation}
By using the analytical solutions in equation \refEquation{eq::internalEnergy}
the exact internal energy for an octant of the spherical shell is obtained as
follows
\begin{equation}
U_{ex} = \dfrac{\pi\,P\,\sigma_y\,r_i^3}{4\,E} \bigg[ ( 1 - \nu )
\dfrac{r_p^3}{r_i^3} -\dfrac{2}{3}( 1 - 2\,\nu )\bigg( 1 - \dfrac{r_p^3}{r_o^3}
+ 3\,\ln\Big(\dfrac{r_p}{r_i} \Big) \bigg) \bigg].
\label{eq::externalEnergy}
\end{equation}
\refFigure{fig::convergenceOfSphericalShell} shows the relative error in energy
with respect to the number of degrees of freedom for studies where the maximum
refinement level is increased from zero to three. To rule out any domain integration
errors, all volumetric integrals were evaluated exactly following~\cite{Kudela2016a}. The curve without any multi-level refinement depicts
exponential convergence in the pre-asymptotic range until the error is dominated
by the kink in the solution field caused by the elastic-plastic interface.
This kink is neither resolved by the boundaries of the finite cells nor is the error controlled by
a refinement in the vicinity of this irregularity in the solution. Therefore, in the asymptotic
range, the error decreases only algebraically when the polynomial degree is further increased. The
kink in the solution field is better approximated by the basis functions as the elements are
refined towards the elastic-plastic interface with hierarchical multi-level refinement. This can be
seen in curves with one and two refinement levels towards the elastic-plastic interface, where the
point at which the convergence levels off to its asymptotic value occurs much later. The curve with
three refinement levels converges exponentially. Here, the error at the interface does not dominate the overall error.

This demonstrates that the combination of the techniques presented
in~\cref{sec:multilevelhp,sec::fcm,sec:elastoplasticity} is able to maintain the expected
convergence behavior for higher-order methods even if neither the physical boundaries of an
artifact nor an evolving elasto-plastic front are explicitly captured by the boundaries of the discretization.

\subsection{Thermo-elasto-plastic bar example} \label{sec:elastoplasticBar}
\begin{figure}[!ht]
\centering
    \includegraphics[]{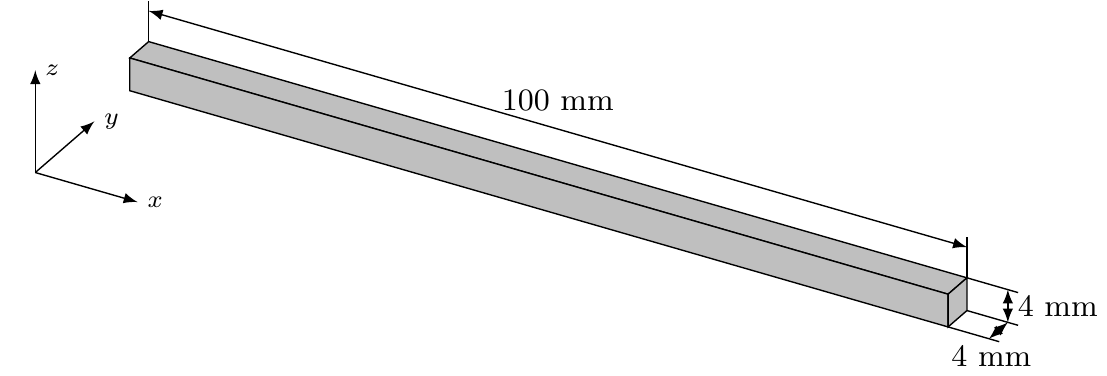}
    \includegraphics[]{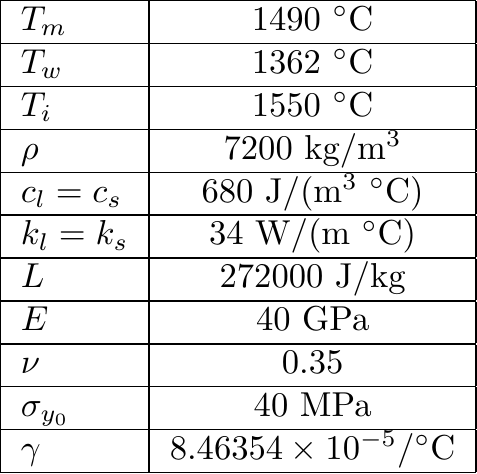}
    \caption{Set up of solidifying bar}%
    \label{fig::solidificationOfABarProblemSetup}%
\end{figure}
In this section, the thermo-mechanical model introduced in
\refSection{sec:coupling} is investigated in the framework of multi-level
\textit{hp}-adaptivity. To this end, an idealized version of the metal
casting process is simulated on a semi-infinite bar. Initially, the bar, which is given in
\refFigure{fig::solidificationOfABarProblemSetup}\,, is liquid with a uniform
temperature, $T_i$. Then, the temperature at the boundary $x=0$ is instantaneously changed to $T_w$ and
kept constant. $T_w$ is lower than the melting
temperature, $T_m$. Therefore, the bar starts to solidify from this boundary onwards.

The analytical solution to the thermal part of this problem is described
in~\cite{Weber1912a} and is known in the literature as \textit{Neumanns's
method}, see e.g.~\cite{Hu1996,Hahn2012}. The position of the liquid-solid
interface is given as
\begin{equation}
X(t) = 2\,\lambda\,\sqrt{\alpha\,t},
\label{eq::interfacePositionBar}
\end{equation}
where $\alpha$ is the diffusivity of the material and $t$ is the time. The
constant $\lambda$ is computed by solving the following nonlinear equation
\begin{equation}
\dfrac{e^{-\lambda^2}}{\text{erf}(\lambda)} +
\dfrac{e^{-\lambda^2}\,(T_m-T_i)}{\text{erfc}(\lambda)\,(T_m-T_w)} =
\dfrac{\lambda\,L\sqrt{\pi}}{c\,(T_m-T_w)},
\label{eq::lambdaEq}
\end{equation}
where $L$ is the latent heat and $c$ is the heat capacity. The analytical
temperature distribution is given in equation
\refEquation{eq::analyticalTemperatureOfBar} for the semi-infinite bar.
\begin{equation}
T(x,t) = 
  \begin{cases}
    T_w + (T_m-T_w) \dfrac{\text{erf}(x/2\sqrt{\alpha
    t})}{\text{erf}(\lambda)} & \quad \text{if} \quad x \leq X(t)\\ \\
    T_i + (T_m-T_i) \dfrac{\text{erfc}(x/2\sqrt{\alpha
    t})}{\text{erfc}(\lambda)}  & \quad \text{if}
    \quad x > X(t)
  \end{cases}
\label{eq::analyticalTemperatureOfBar}
\end{equation}

By using the temperature distribution over the body obtained from the
\textit{Neumanns's method}, Weiner and Boley \cite{Weiner1963} developed an
analytical solution for the thermal stresses for this problem. They assumed that
the body is composed of an elastic perfectly plastic material with a yield
stress, $\sigma_y$, which linearly decreases as the temperature reaches the melting
point, where it becomes zero:
\begin{equation}
\sigma_{y}(T) = 
  \begin{cases}
    \sigma_{y_0}\dfrac{T_m-T}{T_m-T_w} & \quad \text{if} \quad  T_w \leq T \leq
    T_m,\\
    \\
    0 & \quad \text{if} \quad T > T_m.
  \end{cases}
  \label{eq::yieldStressBar}
\end{equation}
Moreover, the following dimensionless quantities are introduced,
\begin{equation}
  m = \frac{(1-\nu)\,\sigma_{y_0}}{\gamma\,E\,(T_m-T_w)}, \quad D =
  \frac{1}{\text{erf}(\lambda)}, \quad \hat{x} = \frac{x}{X(t)},\quad \hat{\sigma} =
  \frac{\sigma\,(1-\nu)}{\gamma\,E\,(T_m-T_w)},
\end{equation}
where $E$, $\nu$ and $\gamma$ are Young's modulus, Poisson's ratio and thermal
expansion coefficient, respectively. The analytical solution of the
dimensionless normal stress components in y- and z-directions are then given as 
\begin{equation}
\hat{\sigma}_{yy}(\hat{x}) = \hat{\sigma}_{zz}(\hat{x}) = 
  \begin{cases}
    m\,[D\,\text{erf}(\lambda\,\hat{x})-1] \hspace{11.2em} \quad \text{if} \quad
    0 \leq \hat{x} < \hat{x}_2\\ \\
    m\,[1-D\,\text{erf}(\lambda\,\hat{x}_1)] + D\,[\text{erf}(\lambda\,
    \hat{x}_1) - \text{erf}(\lambda\,\hat{x})] \\ \quad
    \quad - \dfrac{2}{\sqrt{\pi}}\,D\,(1-m)\,\lambda\,\hat{x}_1\,
    e^{-\lambda^2\hat{x}_1^2}\,
    \text{log}(\dfrac{\hat{x}_1}{\hat{x}}) \hspace{1.3em} \quad \text{if} \quad
    \hat{x}_2 \leq \hat{x} \leq \hat{x}_1 \\ \\
    m\,[1-D\,\text{erf}(\lambda\,\hat{x})] \hspace{11.1em} \quad \text{if}
    \quad \hat{x}_1 < \hat{x} \leq 1 \\ \\
    0 \hspace{18.3em} \quad \text{if} \quad \hat{x} > 1
  \end{cases}
  \label{eq::analyticalStressOfBar}
\end{equation}
where the coordinates $\hat{x}_1$ and $\hat{x}_2$ denote the elastic-plastic
interfaces. According to this solution, the solid material behaves elastic
within the range $[\hat{x}_2,\,\hat{x}_1]$ and plastic in ranges
$[0,\,\hat{x}_2]$ and $[\hat{x}_1,\,1]$. The positions of these interfaces are
obtained by solving the following nonlinear equations
\begin{equation}
\begin{aligned}
2\,(1-m)\,\lambda^2\,\hat{x}_1\,e^{-\lambda^2\hat{x}_1^2}\,(\hat{x}_1-\hat{x}_2)
&=(1+m)\,e^{-\lambda^2\hat{x}_2^2}-(1-m)\,\text{exp}(-\lambda^2\,\hat{x}_1^2)-m(e^{-\lambda^2}+1),
\\
\dfrac{2(1-m)\,\lambda\,\hat{x}_1\,e^{-\lambda^2\hat{x}_1^2}}{\sqrt{\pi}}\,\text{log}(\dfrac{\hat{x}_1}{\hat{x_2}})
&= (1-m)\,\text{erf}(\lambda\,\hat{x}_1)-(1+m)\,\text{erf}(\lambda\,\hat{x}_2)+2\,m\,\text{erf}(\lambda)
\label{eq::ksiEq}
\end{aligned}
\end{equation}

The values that are selected for the material properties of this problem are
provided in \refFigure{fig::solidificationOfABarProblemSetup}\,. With these
values, the constant $\lambda$ is computed as 0.330825295611989, while the
coordinates $\hat{x}_1$ and $\hat{x}_2$ are found to be 0.45487188 and
0.21570439, respectively.

\begin{figure}[!ht]%
  \begin{center}%
    \subfloat[t = 1 second]%
    {
      \includegraphics[width=0.49\textwidth]{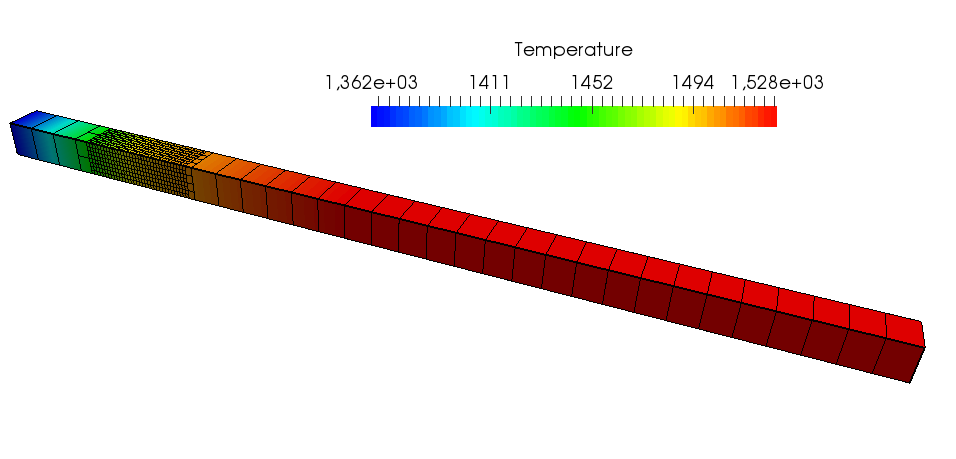}%
      \includegraphics[width=0.49\textwidth]{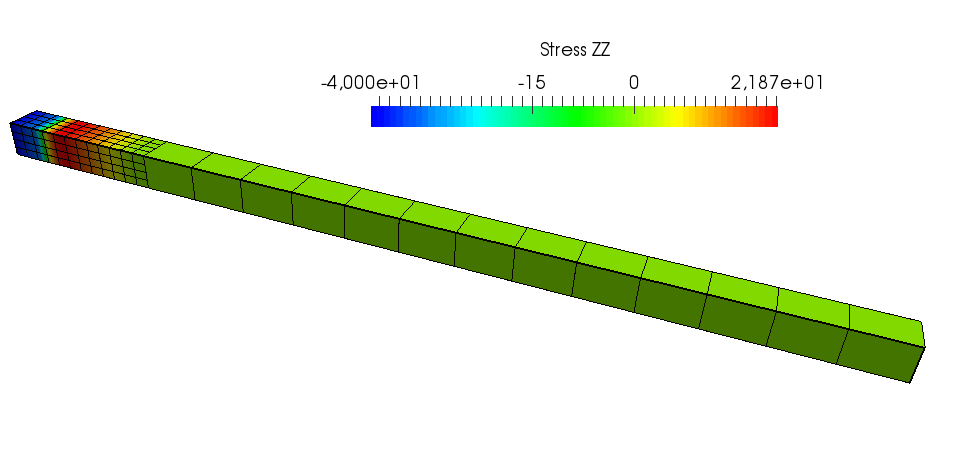}%
    }
    \hfill%
    \subfloat[t = 5.1 seconds]%
    {
      \includegraphics[width=0.49\textwidth]{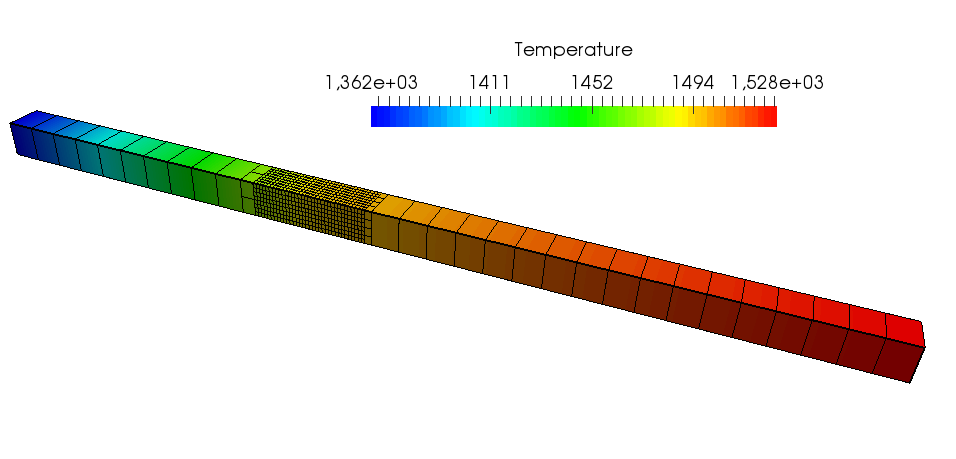}%
      \includegraphics[width=0.49\textwidth]{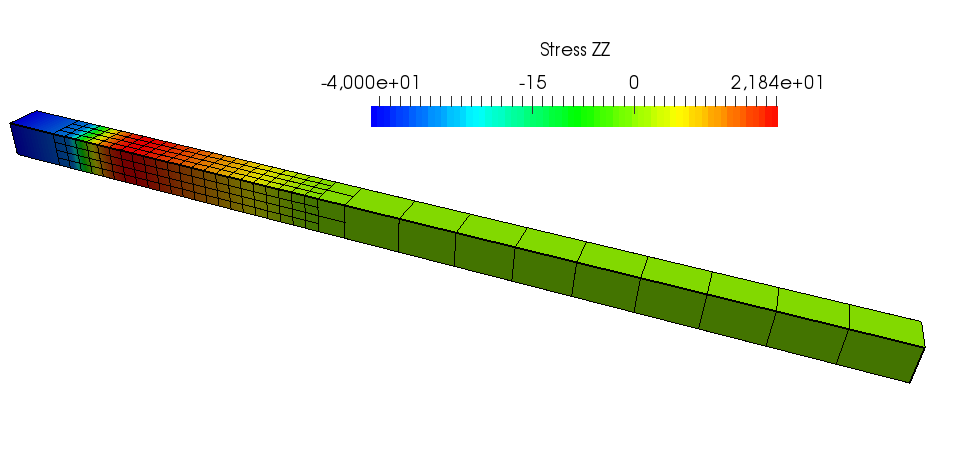}%
    }
    \hfill%
    \subfloat[t = 13.2 seconds]%
    {
      \includegraphics[width=0.49\textwidth]{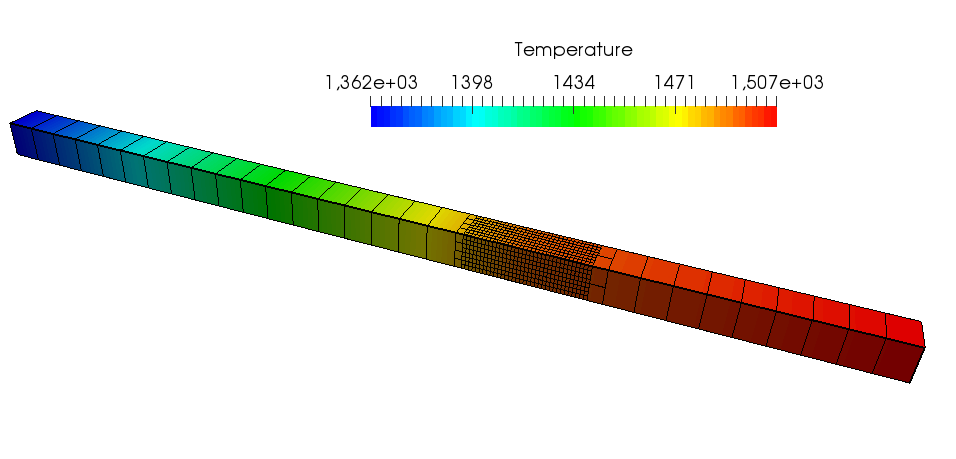}
      \includegraphics[width=0.49\textwidth]{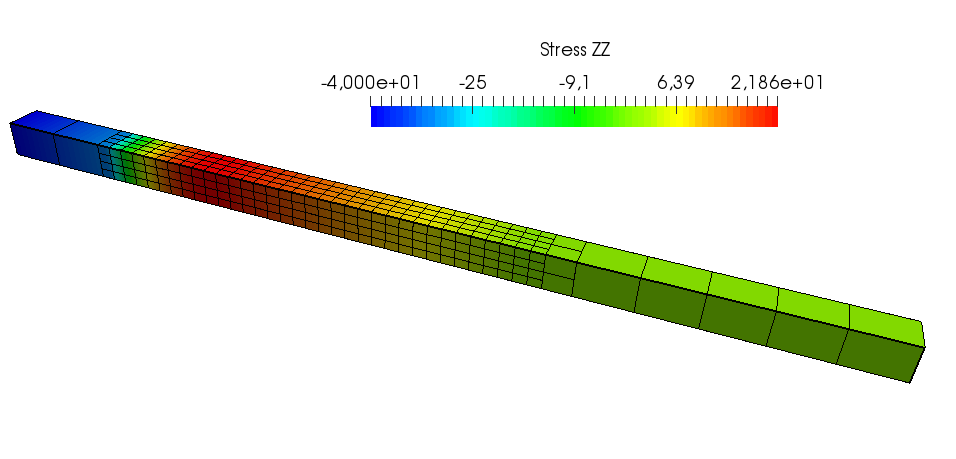}%
    }
    \caption{Temperature and stress distributions at various time states on
    dynamically refined thermal and mechanical meshes, respectively}%
    \label{fig::barResults}%
  \end{center}%
\end{figure}%

The transient thermal problem is solved on a base mesh with 32 hexahedral
elements of order $p=3$ that are distributed in x-direction. Homogeneous
Neumann boundary conditions are applied along y- and z-directions and
Dirichlet boundary conditions are applied in the yz-plane at $x=0$ and $x=100$. In
order to better represent the kink in the temperature field, the elements are
dynamically refined three times towards the solid-liquid interface. The time
step for the backward Euler scheme is chosen to be $\delta t = 0.1$ seconds. 

The simulation for the quasi-static mechanical problem is carried out on a base
mesh with 16 hexahedral elements of order $p=4$ which are dynamically refined
twice towards both the solid-liquid interface and the elastic-plastic
interfaces. During the simulation the liquid part of the bar is treated as a
fictitious domain as explained in \refSection{sec::fcm} by multiplying $E$ with
$\alpha=10^{-8}$. Extended plane strain conditions are applied along y- and
z-directions such that $\varepsilon_{yy}$ and $\varepsilon_{zz}$ are constant.
This is achieved by constraining all degree of freedoms corresponding to
displacements in y- and z-directions to be equal via Lagrange multipliers. As
for the x-direction, the bar is fixed at $x=0$ and left free at $x=100$.

\begin{figure}[!ht]%
  \begin{center}%
    \subfloat[t = 1 second]%
    {
      \includegraphics[]{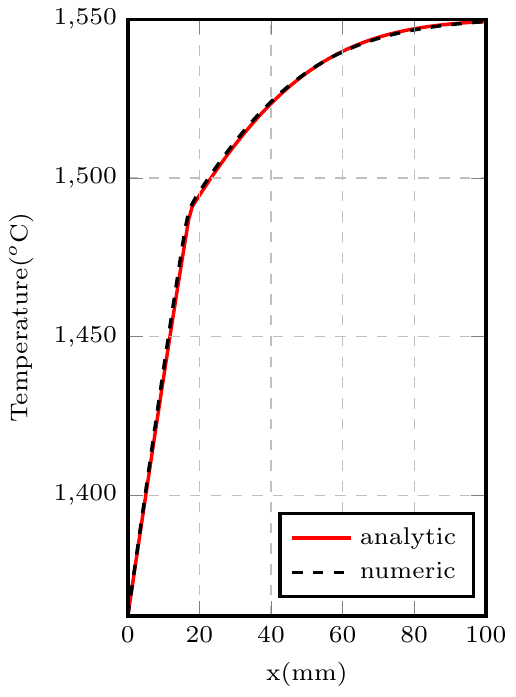}
      }
    \hfill%
    \subfloat[t = 5.1 seconds]%
    {
      \includegraphics[]{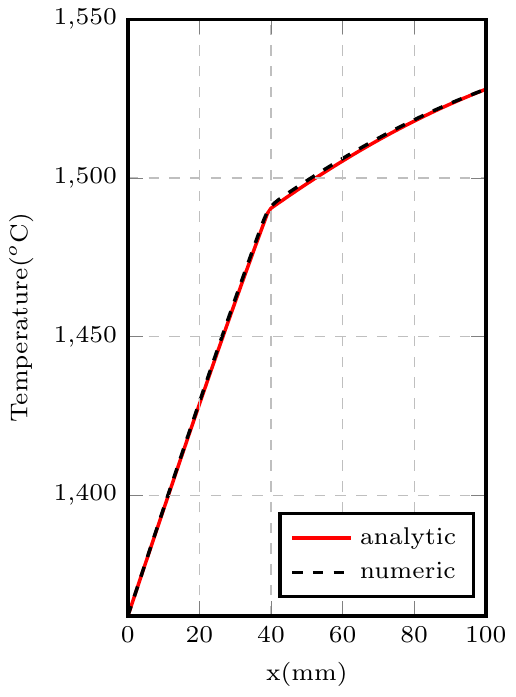}
      }
    \hfill%
    \subfloat[t = 13.2 seconds]%
    {
      \includegraphics[]{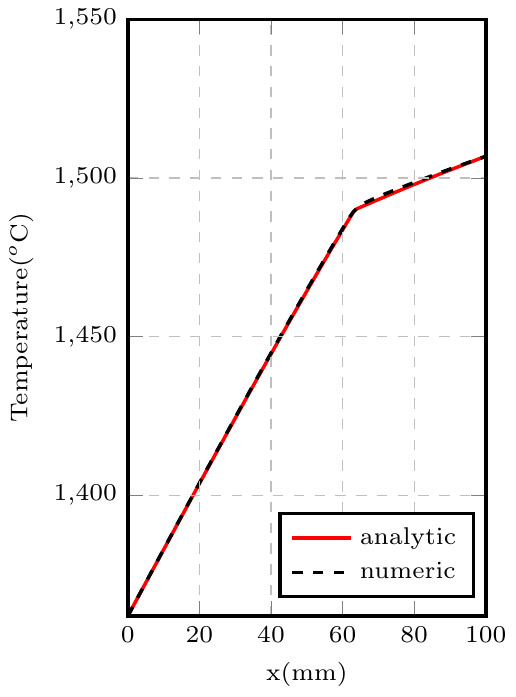}
      }
    \caption{Comparison of the numerical and the analytical temperature
    distribution along x-axis.}%
    \label{fig::temperatureInBar}%
  \end{center}%
\end{figure}%

\begin{figure}[!ht]%
  \begin{center}%
    \subfloat[t = 1 second]%
    {
      \includegraphics[]{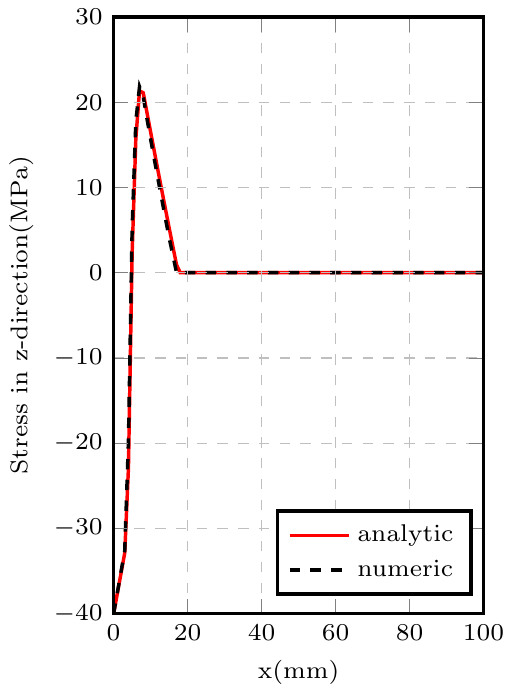}
      }
    \hfill%
    \subfloat[t = 5.1 seconds]%
    {
      \includegraphics[]{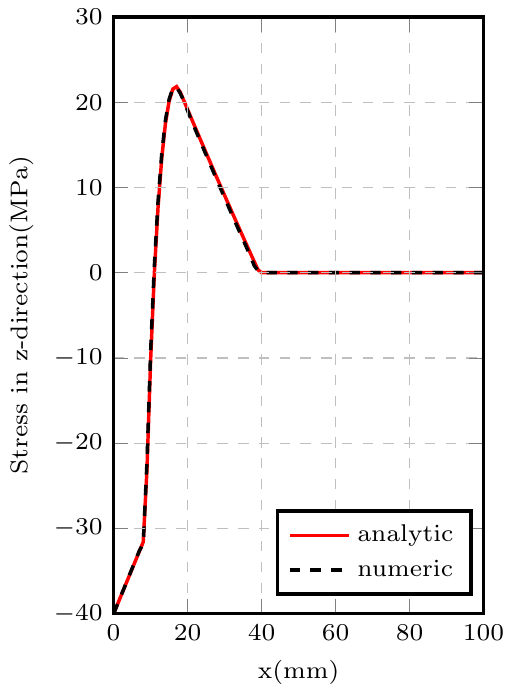}
      }
    \hfill%
    \subfloat[t = 13.2 seconds]%
    {
      \includegraphics[]{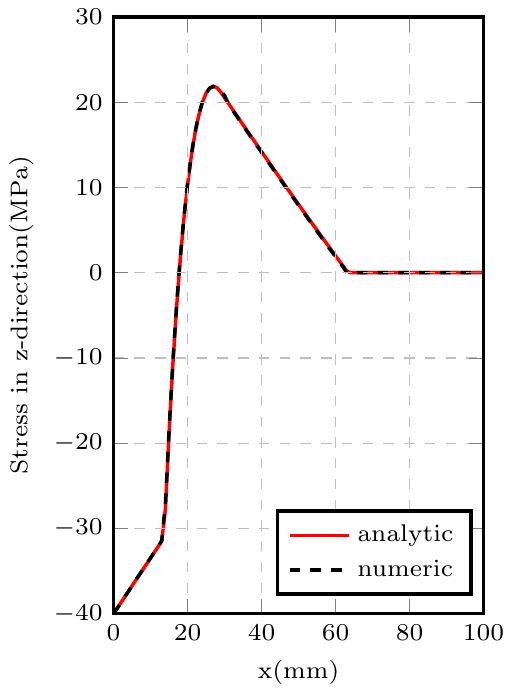}
      }
    \caption{Comparison of the numerical and the analytical stress
    ($\sigma_{zz}$) distribution along x-axis.}%
    \label{fig::stressesInBar}%
  \end{center}%
\end{figure}%

\refFigure{fig::barResults} presents the results of the thermal and the mechanical
problems for the time states $t=\left\{1, 5.1, 13.2\right\}$ seconds. For the thermal problem the evolution of the
temperature field is provided over the dynamically adapted thermal mesh. As for
the mechanical problem the stress component in z-direction ($\sigma_{zz}$) is
given over the dynamically adapted mechanical mesh. The sequence of figures demonstrates how elements are refined and coarsened during the simulation for thermal
and mechanical problems independently.

\refFigure{fig::temperatureInBar} compares the temperature distribution along
the bar to the analytical solution obtained by the \textit{Neumanns's
method}. The kink in the temperature field at melting temperature, which is
caused by the latent heat release during the phase change from liquid to solid,
is well captured in all time states as well as the general temperature profile
along the bar. The stress component in z-direction ($\sigma_{zz}$) along the bar
is depicted in \refFigure{fig::stressesInBar} along with the analytical solution
provided by Weiner and Boley. The kinks at the elastio-plastic interfaces and the
solid-liquid interface are well represented due to refinements. Moreover, it
can be seen that the numerical results closely match their analytical counterparts
along the bar at all time steps of the computation.

\subsection{Applications to metal deposition} \label{sec:metalDepositionExperiment}
\tikzexternaldisable
In this section the performance of the numerical method introduced in
\refSection{sec:numericalMethods} in simulating metal deposition processes is
investigated against the benchmark problem of a single bead laid down
on the top surface of a plate as published in~\cite{Truman2009}. This benchmark
problem is produced by the European Network on Neutron Techniques
Standardization for Structural Integrity (NeT), whose mission is to develop
experimental and numerical techniques and standards for the reliable
characterisation of residual stresses in structural welds. Partners from
industry and academia participated in this benchmark by constructing
an experiment and measuring the residual stresses on the plate \cite{Ohms2009,Ficquet2009,Hofmann2009,Pratihar2009,Turski2009} or
predicting the residual stresses through numerical methods \cite{Shan2009,Ohms2009,Gilles2009,Bate2009,Ficquet2009,Bouchard2009}. Smith
et. al. compares the residual stress measurements and predictions of all these
groups in~\cite{Smith2009}.

\begin{figure}
  \begin{center}
      \includegraphics[]{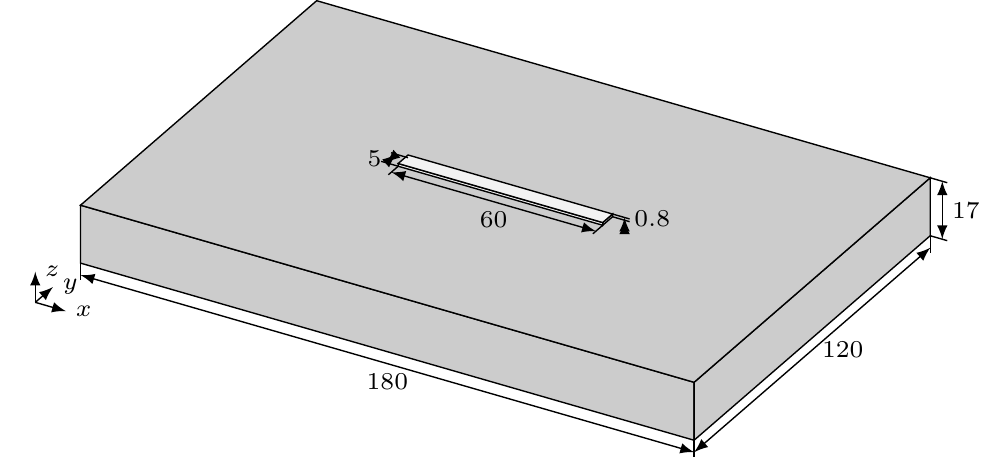}
  \end{center}
  \caption{Setup of welding}
  \label{fig::weldingSetup}
\end{figure}

\refFigure{fig::weldingSetup} illustrates the setup of the benchmark problem
where the plate has dimensions 180 mm $\times$ 120 mm $\times$ 17 mm in
length, width and thickness, respectively. Additionally, the bead that is laid
down on the top surface of the plate has a width of 5 mm, thickness of 0.8 mm
and a length of 60 mm.

According to the benchmark protocol, the plate and the
bead material is AISI Type 316L stainless steel. NeT recommends using slightly
different temperature dependent material properties for the parent and weld
metal. However, in this work the weld metal is not distinguished from the
parent metal. \refTable{tab::weldingSteelProperties} provides the temperature
dependent specific heat, conductivity, thermal expansion coefficient and Young's
modulus. Over the melting temperature of 1400 \degree C these values are
assumed to be constant. In addition to that the density and Poisson's ratio are
assumed to be constant values of 7966 kg m$^{-3}$ and 0.294, respectively and
latent heat is taken as 260 kJ/kg. Finally, the monotonic tensile test
results, which show the true stress corresponding to the percentage of plastic
strain, is given in \refTable{tab::weldingHardeningProperties} by NeT. The
isotropic hardening behaviour of the material is based on this table.

\begin{table}[h!]
  \begin{center}
    \caption{Material properties of AISI Type 316L stainless steel
    \cite{Shan2009}.}
    \label{tab::weldingSteelProperties}
	\begin{tabular}{@{}l cccc@{}}
      \toprule
      \thead{Temperature\\ $[\degree$C]} & \thead{ Specific heat\\
      $[$kJ/kg/\degree C]} & \thead{Conductivity \\ $[$W/m/\degree C]} &
      \thead{Thermal expansion\\ $[\times10^6$ mm/mm/\degree C]} &
      \thead{Young's modulus \\ $[$GPa]}\\
      \midrule
      20  & 0.492& 14.12& 14.56& 195.6\\
      100 & 0.502& 15.26& 15.39& 191.2\\
      200 & 0.514& 16.69& 16.21& 185.4\\
      300 & 0.526& 18.11& 16.86& 179.6\\
      400 & 0.538& 19.54& 17.37& 172.6\\
      500 & 0.550& 20.96& 17.78& 164.5\\
      600 & 0.562& 22.38& 18.12& 155.0\\
      700 & 0.575& 23.81& 18.43& 144.1\\
      800 & 0.587& 25.23& 18.72& 131.4\\
      900 & 0.599& 26.66& 18.99& 116.8\\
      1000& 0.611& 28.08& 19.27& 100.0\\
      1100& 0.623& 29.50& 19.53& 80.0\\
      1200& 0.635& 30.93& 19.79& 57.0\\
      1300& 0.647& 32.35& 20.02& 30.0\\
      1400& 0.659& 33.78& 20.21& 2.0\\
      \bottomrule
    \end{tabular}
  \end{center}
\end{table}

\begin{table}[h!]
  \begin{center}
    \caption{True stress values [MPa] of AISI Type 316L stainless steel at
    various true plastic strain \cite{Ficquet2009}.}
    \label{tab::weldingHardeningProperties}
	\begin{tabular}{@{}l cccccccccc@{}}
      \toprule
      \thead{Temperature $[\degree$C]} & \thead{0\%} & \thead{0.2\%} &
      \thead{1\%} & \thead{2\%} & \thead{5\%} & \thead{10\%} &
      \thead{20\%} & \thead{30\%} & \thead{40\%}\\
      \midrule
      23  & 210 & 238  & 292& 325& 393& 494& 648& 775& 880\\
      275 & 150 & 173.7& 217& 249& 325& 424& 544& 575& \\
      550 & 112 & 142.3& 178& 211& 286& 380& 480& 500& \\
      750 & 95  & 114.7& 147& 167& 195& 216& 231& 236& \\
      800 & 88  & 112  & 120& 129& 150& 169& 183&    & \\
      900 & 69  & 70   & 71 & 73 & 76 & 81 &    &    & \\
      1100& 22.4&      &    &    &    &    &    &    & \\
      1400& 2.7 &      &    &    &    &    &    &    & \\
      \bottomrule
    \end{tabular}
  \end{center}
\end{table}

The heat input from the welding torch, which travels with the speed of 2.27
mm$/$s, is defined by NeT as 633 J$/$mm. The moving heat source is based on the
Goldak's \cite{Goldaknewfiniteelement1984} ellipsoidal model, which is widely
used in welding simulations. In this model the body load is given as
\begin{equation}
q(x,y,z) =
\dfrac{6\,\sqrt{3}\,Q}{\pi\,\sqrt{\pi}\,r_x\,r_y\,r_z}\,e^{-3\Big(\dfrac{x-x_c}{r_x}\Big)^2}\,e^{-3\Big(\dfrac{y-y_c}{r_y}\Big)^2}\,e^{-3\Big(\dfrac{z-z_c}{r_z}\Big)^2},
\label{eq::goldakHeatSource}
\end{equation}
where $x_c$, $y_c$ and $z_c$ define the position of the welding torch and $r_x$,
$r_y$ and $r_z$ are the semi-axes of the ellipsoid. Shan et.
al.\cite{Shan2009} uses the weld bead profile measurements to decide the
semi-axes of the ellipsoid, which is followed in this work. According to these
measurements $r_x$, $r_y$ and $r_z$ are assumed to be 1.9 mm, 3.2 mm and 2.8 mm,
respectively. Moreover, in equation \refEquation{eq::goldakHeatSource} $Q$
is the power of the heat source, which is computed from the heat input and the
welding speed as 1437 Watt.

The heat loss from the system is modeled with radiation and convection
boundary conditions on the surface of the weld bead and the plate. During the
process, the surface is updated as new weld bead is deposited. The convection
coefficient and the emissivity are taken as 10 W$/$m$^2$K and 0.75, respectively
and they are assumed to be independent of the temperature. Additionally, the
ambient temperature is also assumed to be constant at 20 \degree C. 

\begin{figure}
  \begin{center}     
      \includegraphics[]{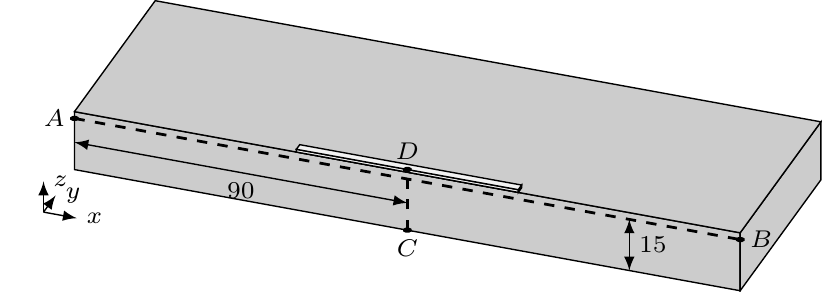}
  \end{center}
  \caption{Setup of welding}
  \label{fig::weldingSetupHalf}
\end{figure}

Due to symmetry, only half of the plate and the weld bead is considered as shown
in \refFigure{fig::weldingSetupHalf} for the numerical model with the appropriate
symmetry boundary conditions. The plate and the weld bead are embedded in a
Cartesian grid of size $30\times12\times7$ elements with polynomial order 3 for
both the thermal and mechanical problems. The dimensions of the computational
domain are 180 mm, 60 mm and 18.6 mm in length, width and thickness, respectively. It should be noted
that the thickness of the computational domain is greater than the total
thickness of the plate and the weld bead, which is 17.8 mm. The first six rows of
elements in thickness are used to discretize the plate, while the weld bead is
embedded in the last row of elements. 

\begin{figure}%
  \begin{center}%
    \subfloat[Mesh of the thermal problem]%
    {
      \includegraphics[width=0.45\textwidth]{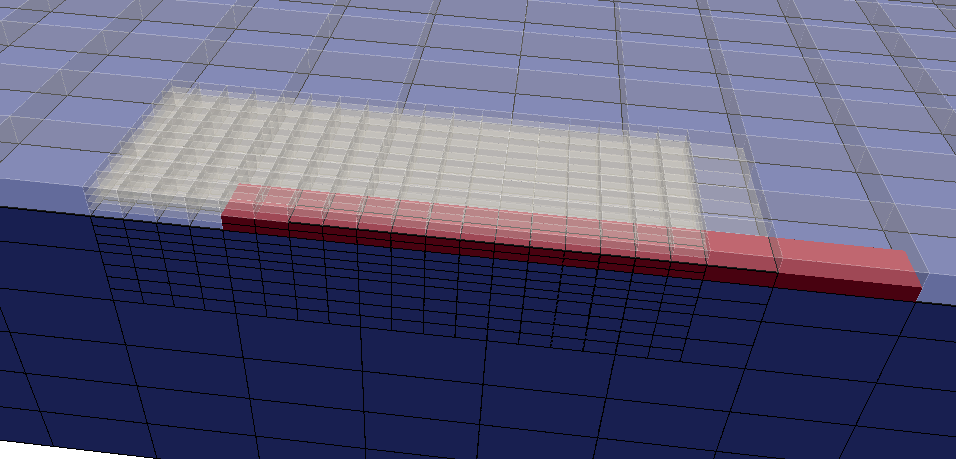}
    }
    \hfill%
    \subfloat[Mesh of the mechanical problem]%
    {
      \includegraphics[width=0.45\textwidth]{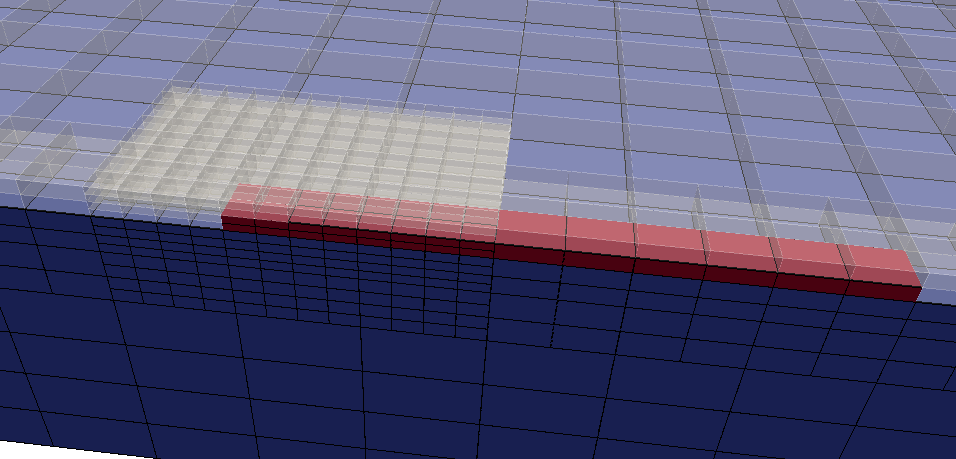}
    }
    \caption{Meshes used for numerical simulation with material states, where
    red, blue and white (hollow) denote weld bead, plate and air, respectively.
    }%
    \label{fig::weldingFCM}%
  \end{center}%
\end{figure}%

\begin{figure}%
  \begin{center}%
    \subfloat[t = 1 second]%
    {
      \includegraphics[width=0.45\textwidth]{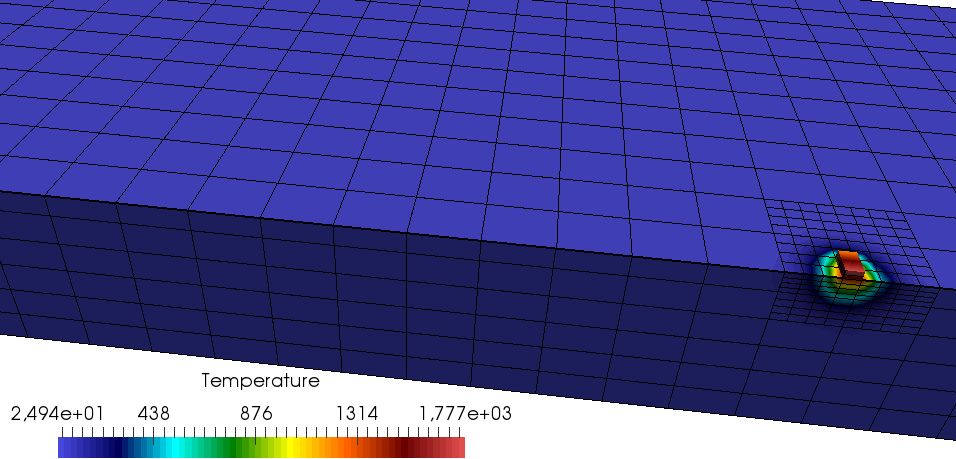}
      \includegraphics[width=0.45\textwidth]{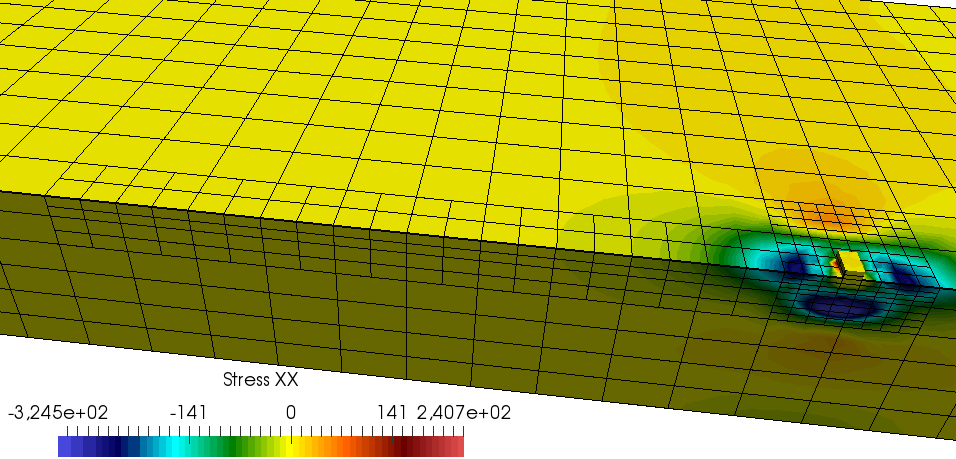}
    }
    \hfill%
    \subfloat[t = 5.1 seconds]%
    {
      \includegraphics[width=0.45\textwidth]{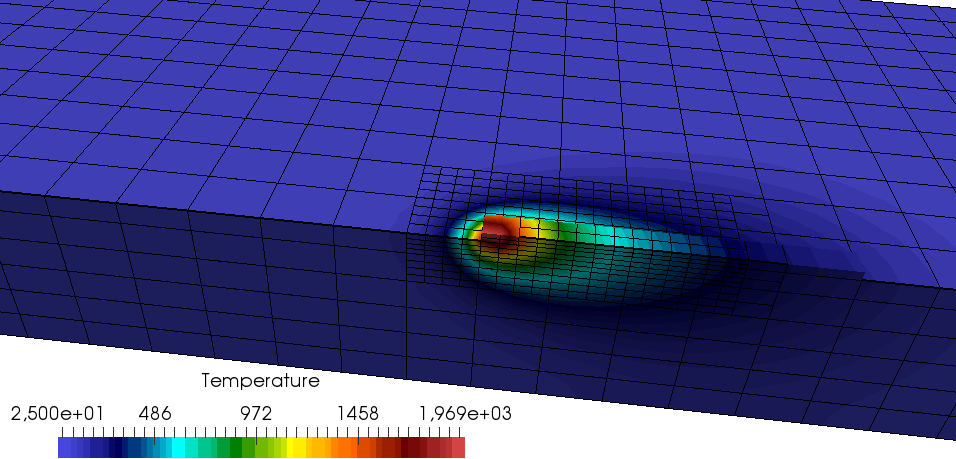}
      \includegraphics[width=0.45\textwidth]{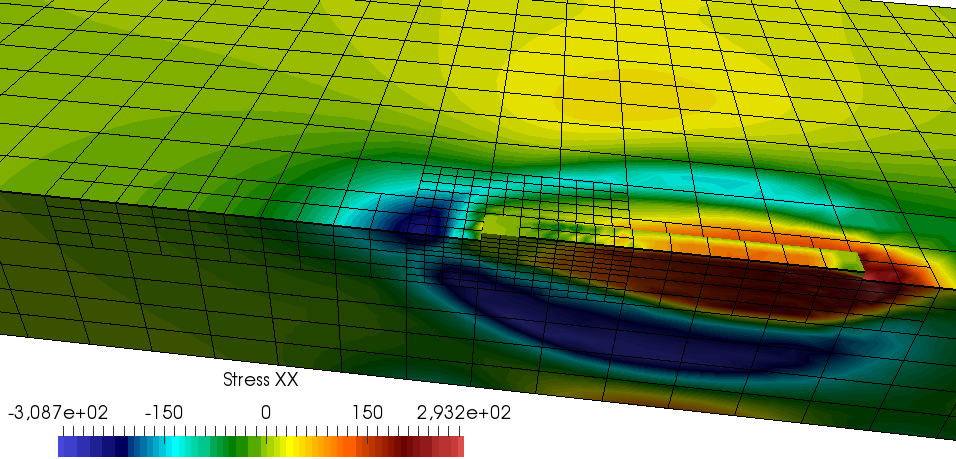}
    }
    \hfill%
    \subfloat[t = 13.2 seconds]%
    {
      \includegraphics[width=0.45\textwidth]{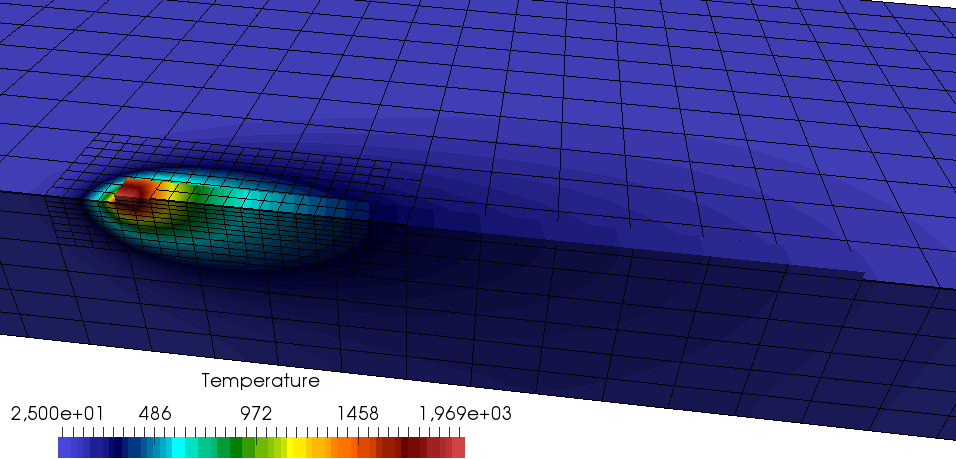}
      \includegraphics[width=0.45\textwidth]{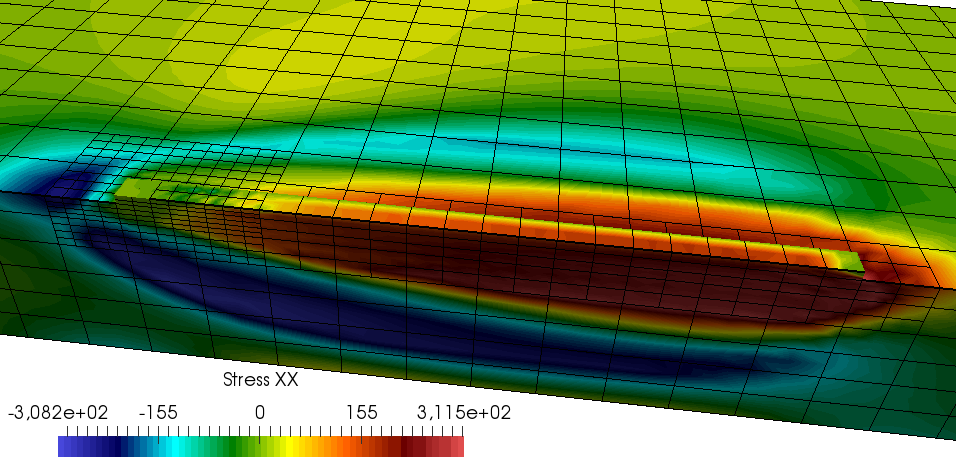}
    }
    \caption{Evolution of the temperature and longitudinal stress longitudinal
    stress ($\sigma_{xx}$) distributions during the welding process on dynamic
    meshes at various time states.}%
    \label{fig::weldingResults}%
  \end{center}%
\end{figure}

\refFigure{fig::weldingFCM} shows a snapshot of the dynamic meshes of the
thermal and the mechanical problems with material states, where red, blue and
white (hollow) denote weld bead, plate and air, respectively. It should be noted
that only the mechanical mesh is refined once along the welding line before
deposition starts, which is not coarsened afterwards. Whereas, both meshes are
dynamically refined twice towards the weld torch when material is deposited.
The weld torch is set to travel 1.5 mm in each time step. Therefore, the newly
deposited domain near the weld torch can be captured by the sub-element
boundaries and the multi-level grid that keeps track of the physical domain. As the
weld torch moves further away, these sub-elements are coarsened while the
multi-level grid is used to remember the deposited domain as explained in
\refSection{sec:coupling}\,. \refFigure{fig::weldingResults} depicts the dynamic
meshes with temperature and longitudinal stress ($\sigma_{xx}$) distributions
during the deposition period. 

\begin{figure}%
  \begin{center}%
    \subfloat[Longitudinal stresses ($\sigma_{xx}$)]%
    {
      \includegraphics[width=0.48\textwidth]{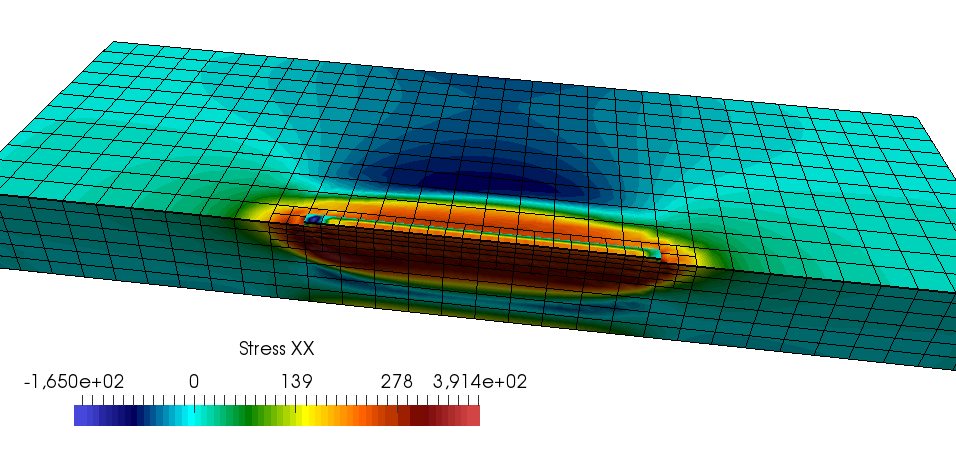}
      }
    \hfill%
    \subfloat[Transverse stresses ($\sigma_{yy}$)]%
    {
      \includegraphics[width=0.48\textwidth]{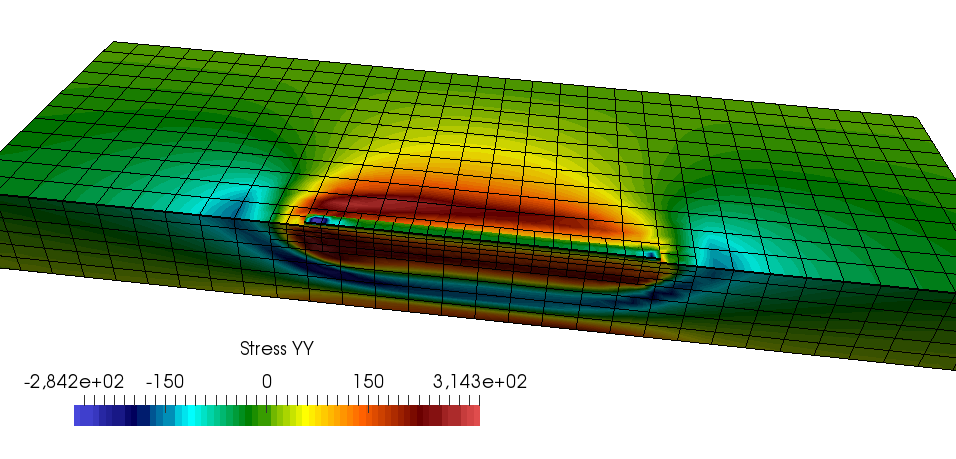}
      }
    \caption{Residual stress distributions after cooling down to room
    temperature.}%
    \label{fig::weldingResidual}%
  \end{center}%
\end{figure}%

\begin{figure}%
  \begin{center}%
    \subfloat[Longitudinal stresses ($\sigma_{xx}$)]%
    {
      \includegraphics[]{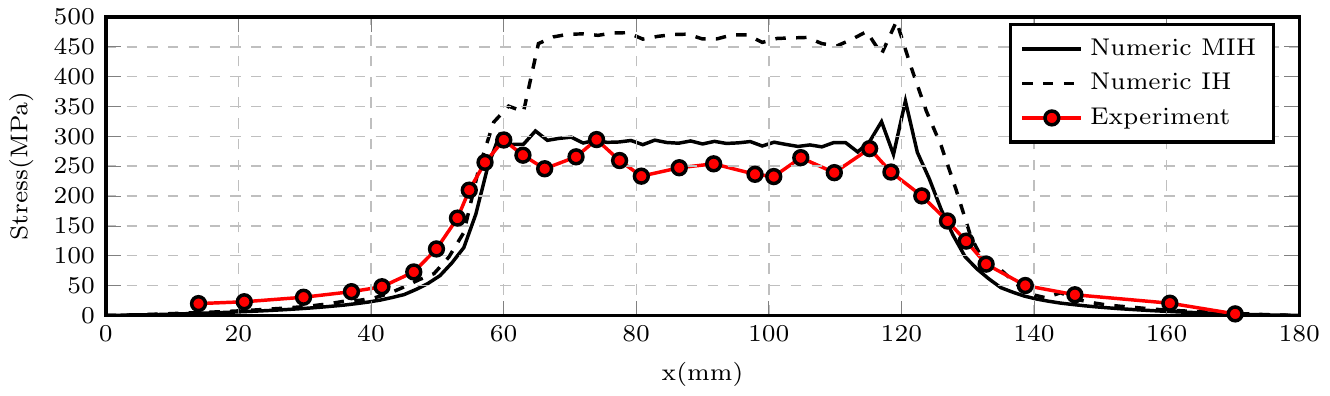}
      }
    \hfill%
    \subfloat[Transverse stresses ($\sigma_{yy}$)]%
    {
      \includegraphics[]{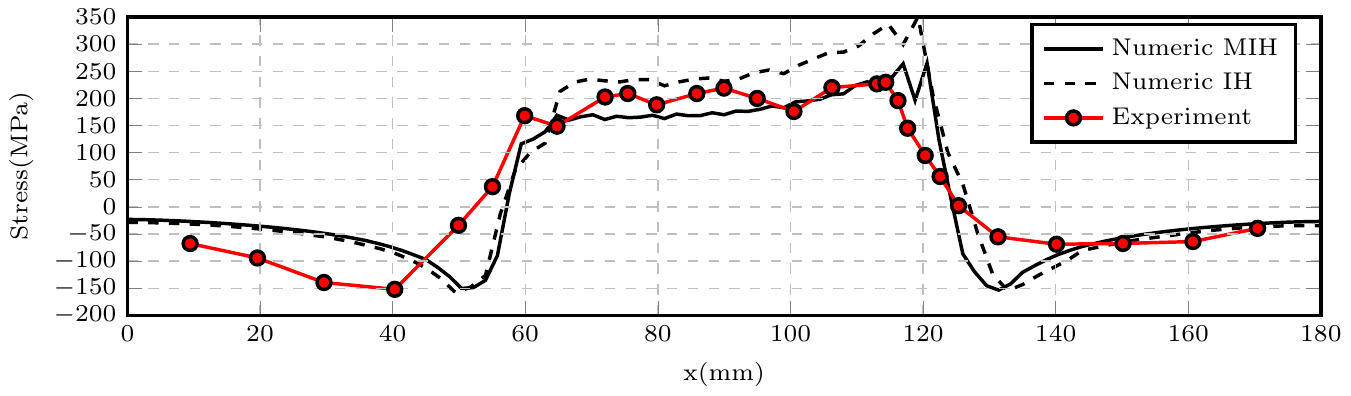}
      }
    \caption{Numerical and experimental (\cite{Gilles2009}) results over line
    AB (IH: Isotropic hardening, MIH: Modified isotropic hardening)}%
    \label{fig::plotsOveLineAB}%
  \end{center}%
\end{figure}%

\begin{figure}%
  \begin{center}%
    \subfloat[Longitudinal stresses ($\sigma_{xx}$)]%
    {
      \includegraphics[]{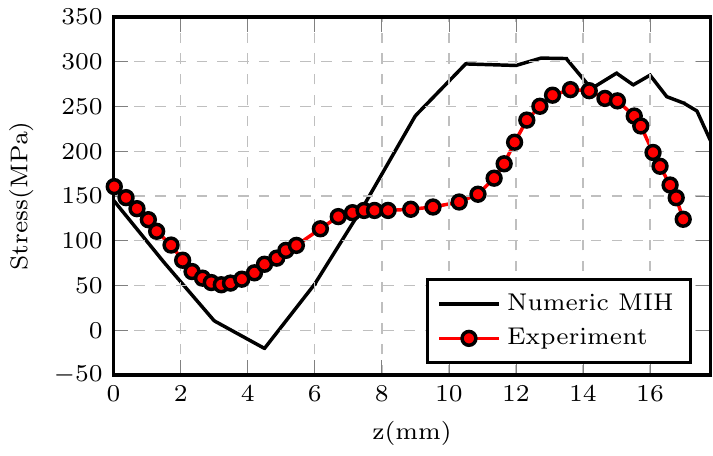}
      }
    \hfill%
    \subfloat[Transverse stresses ($\sigma_{yy}$)]%
    {
      \includegraphics[]{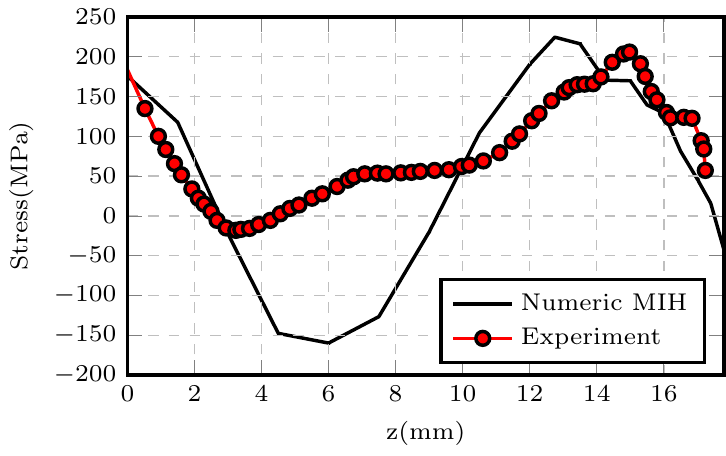}
      }
    \caption{Numerical and experimental (\cite{Gilles2009}) results over line
    CD (MIH: Modified isotropic hardening)}%
    \label{fig::plotsOveLineCD}%
  \end{center}%
\end{figure}%

Residual stress distributions after the parts are cooled down to the room
temperature are given in \refFigure{fig::weldingResidual} over the physical
domain. Gilles et. al. \cite{Gilles2009} measures the longitudinal and
transverse residual stresses with a neutron diffraction technique over the lines AB
and CD. These lines can be seen in \refFigure{fig::weldingSetupHalf}\,.
\refFigure{fig::plotsOveLineAB} and \refFigure{fig::plotsOveLineCD} compare the
numerically predicted residual stresses to the measurements from experiment by Gilles et. al. along the
lines AB and CD, respectively. It can be seen from these figures that the
measured and the numerically computed stress profiles are in good agreement.
However, utilizing the complete isotropic hardening values in \refTable{tab::weldingHardeningProperties} results
in higher residual stress especially for the longitudinal stress component,
which is also reported by Gilles et. al. \cite{Gilles2009}. Their explanation
for effect is that the isotropic hardening law which does not consider stress
relaxation due to viscoplastic effects or hardening recovery is not good enough
to simulate cyclic loads which occur during the welding process. Therefore, the
isotropic hardening model is modified in this work such that over one percent of
plastic strain the material is assumed to behave elastic perfectly plastic. The
results with this modified model match the experimental values as illustrated in
\refFigure{fig::plotsOveLineAB}.

\section{Summary and conclusions} \label{sec:conclusions}
The article at hand evaluates if high-order finite elements may be used in metal deposition modeling. To this end, the finite cell method, an embedded domain method of high order, was used. It was demonstrated that this method, combined with the multi-level $hp$-method, leads to high-order convergence rates even if neither, the physical boundaries nor the plastic front is resolved by the boundaries of the computational mesh. Therefore, high-order embedded domain modeling is a valid computational option. Very accurate results were also achieved in the thermo-elasto-plastic setting. Therein, the thermal and the elasto-plastic field variables were computed on their own high-order discetizations and coupled by means of a classical staggered scheme. In all of the computations the state variables were decoupled from the finite cells in an FCM sense by managing them on another computational grid. This allows for their separate coarsening and refinement and a simulation where material can be submerged into the computation on a sub-element level i.e. the addition of material is not equal to the addition of finite cells. The methods were thoroughly verified in semi-analytical benchmarks and then applied to an experimental benchmark of metal deposition. A good agreement with the experimentally measured residual strains/stresses was found. 

While the overall computational methodology leads to very accurate results open research questions remain. For example, the presented de-refinements in the mechanical part of the simulation are clearly a valid option to accurately compute the stress states in the vicinity of a plastic front. However, the stress state in a manufactured part does not only consist of one plastic front and rather smooth stress states to the side of it. 
Instead, very complex, local stress states may be present. A straightforward coarsening procedure of the grid in these regions will cause smoothing effects for variables which are not diffusive by nature. 
It is beyond the scope of this paper to study the influence of strong, local coarsening in complex stress states and their effect in more involved computations. The investigation of this challenge remains open. 
Since the local plastic strains are not connected to the finite cells but stored in an extra grid which acts like a computational database, one option would be not to coarsen them at all. 
The computational grid, carrying the degrees of freedom could be coarsened independently. Such a procedure would allow for a constant number of degrees of freedom throughout the computation 
while detailed information on the local plastic strains would still be available if needed even if the finite cells were coarsened. 

Next, it is worthwile to investigate the wall-clock time advantages of the computational approach w.r.t. conventional modeling techniques in  detail. This is not easy to some extent because 
a fair comparison requires an efficient implementation of both techniques which is not available at present. Clearly, the modeling flexibility gained by the presented embedded domain method of high order 
using transient refinements does introduce a certain implementational as well as a computational overhead as compared to a classical fixed-grid approaches. 
The  break-even point of refining and de-refining regarding wall-clock time further depends on the scale of the computations. However, to refine everywhere to the finest level required locally is also not an option.
 Thus, refinements and de-refinements must be used for efficient computations. How this is possible was investigated in the article at hand, but its comparison to other refinement techniques remains an open question.
  The current computations hint that the type of microscopic computations in the metal depositing process presented in the last example pay off in a multi-layer process whereby the refinement and de-refinement is not 
  carried out at each time increment. 

In summary, the article presents the first verified and validated steps of a novel computational methodology for the analysis of metal deposition modeling in an embedded domain sense using locally refined, 
transient discretizations of $hp$-type. As such, it demonstrates that this computational methodology is a valid option for the analysis of metal deposition.



\section*{Acknowledgements} 
The authors gratefully acknowledge the financial support of the German Research Foundation (DFG) under Grant RA 624/27-2.

\clearpage
\section*{References}
\bibliographystyle{elsarticle-num}
\bibliography{SimulationInAppliedMechanicsatCiE}

\end{document}